

  \newcount\fontset
  \fontset=1
  \def\dualfont#1#2#3{\font#1=\ifnum\fontset=1 #2\else#3\fi}

  \dualfont\bbfive{bbm5}{cmbx5}
  \dualfont\bbseven{bbm7}{cmbx7}
  \dualfont\bbten{bbm10}{cmbx10}
  \font \eightbf = cmbx8
  \font \eighti = cmmi8 \skewchar \eighti = '177
  \font \eightit = cmti8
  \font \eightrm = cmr8
  \font \eightsl = cmsl8
  \font \eightsy = cmsy8 \skewchar \eightsy = '60
  \font \eighttt = cmtt8 \hyphenchar\eighttt = -1
  \font \msbm = msbm10
  \font \sixbf = cmbx6
  \font \sixi = cmmi6 \skewchar \sixi = '177
  \font \sixrm = cmr6
  \font \sixsy = cmsy6 \skewchar \sixsy = '60
  \font \tensc = cmcsc10
  
  \font \titlefont = cmr7 scaled \magstep4
  \scriptfont \bffam = \bbseven
  \scriptscriptfont \bffam = \bbfive
  \textfont \bffam = \bbten

  \font\rs=rsfs10 

  \newskip \ttglue

  \def \eightpoint {\def \rm {\fam0 \eightrm }%
  \textfont0 = \eightrm
  \scriptfont0 = \sixrm \scriptscriptfont0 = \fiverm
  \textfont1 = \eighti
  \scriptfont1 = \sixi \scriptscriptfont1 = \fivei
  \textfont2 = \eightsy
  \scriptfont2 = \sixsy \scriptscriptfont2 = \fivesy
  \textfont3 = \tenex
  \scriptfont3 = \tenex \scriptscriptfont3 = \tenex
  \def \it {\fam \itfam \eightit }%
  \textfont \itfam = \eightit
  \def \sl {\fam \slfam \eightsl }%
  \textfont \slfam = \eightsl
  \def \bf {\fam \bffam \eightbf }%
  \textfont \bffam = \eightbf
  \scriptfont \bffam = \sixbf
  \scriptscriptfont \bffam = \fivebf
  \def \tt {\fam \ttfam \eighttt }%
  \textfont \ttfam = \eighttt
  \tt \ttglue = .5em plus.25em minus.15em
  \normalbaselineskip = 9pt
  \def \MF {{\manual opqr}\-{\manual stuq}}%
  \let \sc = \sixrm
  \let \big = \eightbig
  \setbox \strutbox = \hbox {\vrule height7pt depth2pt width0pt}%
  \normalbaselines \rm }


  \def \Headlines #1#2{\nopagenumbers
    \advance \voffset by 2\baselineskip
    \advance \vsize by -\voffset
    \headline {\ifnum \pageno = 1 \hfil
    \else \ifodd \pageno \tensc \hfil \lcase {#1} \hfil \folio
    \else \tensc \folio \hfil \lcase {#2} \hfil
    \fi \fi }}

  \def \Title #1{\vbox{\baselineskip 20pt \titlefont \noindent #1}}

  \def \Date #1 {\footnote {}{\eightit Date: #1.}}

  \def \Authors #1{\bigskip \bigskip \noindent #1}

  \long \def \Addresses #1{\begingroup \eightpoint \parindent0pt
\medskip #1\par \par \endgroup }

  \long \def \Abstract #1{\begingroup \eightpoint
  \bigskip \bigskip \noindent
  {\sc ABSTRACT.} #1\par \par \endgroup }


  \def \lcase #1{\edef \auxvar {\lowercase {#1}}\auxvar }

  \def \goodbreak {\vskip0pt plus.1\vsize \penalty -250 \vskip0pt
plus-.1\vsize }

  \newcount \secno \secno = 0
  \newcount \stno

  \def \seqnumbering {\global \advance \stno by 1
    \number \secno .\number \stno }

  \def \label #1{\def\localvariable {\number \secno
    \ifnum \number \stno = 0\else .\number \stno \fi }\global \edef
    #1{\localvariable }}

  \def\section #1{\global\def\SectionName{#1}\stno = 0 \global
\advance \secno by 1 \bigskip \bigskip \goodbreak \noindent {\bf
\number \secno .\enspace #1.}\medskip \noindent \ignorespaces}

  \long \def \sysstate #1#2#3{\medbreak \noindent {\bf \seqnumbering
.\enspace #1.\enspace }{#2#3\vskip 0pt}\medbreak }
  \def \state #1 #2\par {\sysstate {#1}{\sl }{#2}}
  \def \definition #1\par {\sysstate {Definition}{\rm }{#1}}


  \def \proof {\medbreak \noindent {\it Proof.\enspace }}
  \def \proofend {\ifmmode \eqno \square \else \hfill \square
\looseness = -1 \medbreak \fi }

  \def \$#1{#1 $$$$ #1}
  \def\=#1{\buildrel (#1) \over =}

  \def\Item #1{\smallskip \item {#1}}
  \newcount \zitemno \zitemno = 0
  \def\izitem {\zitemno = 0}
  \def\zitem {\global \advance \zitemno by 1 \Item {{\rm(\romannumeral
\zitemno)}}}

  \newcount \footno \footno = 1
  \newcount \halffootno \footno = 1
  \def\footcntr {\global \advance \footno by 1
  \halffootno =\footno
  \divide \halffootno by 2
  $^{\number\halffootno}$}


  \def \({\left (\vrule height 9pt width 0pt}
  \def \){\right )}
  \def \[{\left \Vert }
  \def \]{\right \Vert }
  \def \*{\otimes }
  \def \+{\oplus }
  \def \:{\colon }
  \def \<{\left \langle }
  \def \>{\right \rangle }
  \def \text #1{\hbox {\rm #1}}
  \def \curly#1{\hbox{\rs #1\/}}
  \def \ds{\displaystyle}
  \def \and {\hbox {,\quad and \quad }}
  
  \def \calcat #1{\,{\vrule height8pt depth4pt}_{\,#1}}
  \def \labelarrow#1{\ {\buildrel #1 \over \longrightarrow}\ }
  
  \def \crossproduct {{\hbox {\msbm o}}}
  
  \def \for #1{,\quad \forall\,#1}
  \def \inv {^{-1}}
  \def \pmatrix #1{\left [\matrix {#1}\right ]}
  \def \square {\hbox {$\sqcap \!\!\!\!\sqcup $}}
  \def \stress #1{{\it #1}\/}

  \def \|{\Vert }
  \def \inv {^{-1}}


  \newcount \bibno \bibno =0
  \def \newbib #1{\global \advance \bibno by 1 \edef #1{\number
    \bibno}}
  \def\cite #1{{\rm [\bf #1\rm ]}}
  \def\scite #1#2{\cite{#1{\rm \hskip 0.7pt:\hskip 2pt #2}}}
  \def\lcite #1{(#1)}
  \def\fcite #1#2{\lcite{#1}}
  \def\bibitem#1#2#3#4{\smallskip \item {[#1]} #2, ``#3'', #4.}

  \def \references {
    \begingroup
    \bigskip \bigskip \goodbreak
    \eightpoint
    \centerline {\tensc References}
    \nobreak \medskip \frenchspacing }


  \def\map#1{\mathrel{\buildrel #1 \over \longrightarrow}}
  \def\b{\beta}
  \def\N{{\bf N}}
  \def\Ctr{{\cal Z}}
  \def\R{{\curly R}}
  \def\Reals{{\bf R}}
  \def\E{{\cal E}}
  \def\Toep{{\cal T}(\R,\E)}
  \def\Alg{C^*(\R,\E)}
  \def\inclow#1{\underline{\vrule depth 1.5pt width 0pt#1}}
  \def\inc#1{\underline{#1}}
  \def\K{{\cal K}}
  \def\Kt{{\hat \K}}
  \def\Kb{{\check\K}}
  \def\ind{\mathop {\rm ind}}
  \def\et{{\hat e}}
  \def\eb{{\check e}}
  \def\tphi{{\hat \varphi}}
  
  \def\d{\, d}
  \def\CP{A\mathop{\crossproduct_{\a,\Tr}} {\bf N}}
  \def\Tr{{\cal L}}
  \def\a{\alpha}
  \def\hb{_{_{H,\b}}}
  \def\ev{c}
  \def\frac#1#2{{#1 \over #2}}


  \def\Bo{Bo}
  \def\Baladi{Ba}
  \def\BS{BS}
  \def\BR{BR}
  \def\CLT{CLT}
  \def\Craizer{C}
  \def\endo{E1}
  \def\tower{E2}
  \def\Fan{F}
  \def\FanJiangOne{FJ1}
  \def\FanJiangTwo{FJ2}
  \def\FL{FL}
  \def\H{H}
  \def\Keller{K}
  \def\Lone{L1}
  \def\Ltwo{L2}
  \def\Lthree{L3}
  \def\LT{LT}
  \def\M{M}
  \def\Ped{P}
  \def\RenaultThesis{Re1}
  \def\RenaultAF{Re2}
  \def\RenaultRadon{Re3}
  \def\Rone{Ru1}
  \def\Rdois{Ru2}
  \def\Rtres{Ru3}
  \def\Wone{W1}
  \def\Wtwo{W2}
  \def\Watatani{Wa}
  \def\Y{Y}

  \def\titletext{C*-Algebras, Approximately Proper Equivalence
  Relations, and
Thermodynamic Formalism}

  \Headlines
  {Equivalence Relations and Thermodynamic Formalism}
  {R.~Exel and A.~Lopes}

  \Title{\titletext}

  \Date{7 Jun 2002}

  \Authors
  {R.~Exel\footnote{*}{\eightrm Partially supported by CNPq.}
  and
  A.~Lopes\footnote{**}{\eightrm Partially supported by CNPq.}}

  \Addresses
  {Departamento de Matem\'atica,
  Universidade Federal de Santa Catarina,
  Florian\'opolis,
  Brazil
  (exel@mtm.ufsc.br),
  \par
  Instituto de Matem\'atica, 
  Universidade Federal do Rio Grande do Sul,
  Porto Alegre, Brazil
  (alopes@mat.ufrgs.br).}

  \Abstract {We introduce a non-commutative generalization of the
notion of (approximately proper) equivalence relation and propose the
construction of a ``quotient space''.
  We then consider certain one-parameter groups of automorphisms of
the resulting C*-algebra and prove the existence of KMS states at
every temperature.  In a model originating from Thermodynamics we
prove that these states are unique as well.  
  We also show a relationship between maximizing measures (the
analogue of the Aubry-Mather measures for expanding maps) and ground
states.  In the last section we explore an interesting example of
phase transitions.}

  \begingroup
  \long\def\ll#1:#2:#3;{\noindent\hskip 0.1\hsize \hbox to 0.8\hsize{
    \hbox to 22pt{\hfil #1.\ \ }#2 \dotfill \hbox to 13pt{\hfil #3}}\par}
  \long\def\ni#1:#2;{\noindent\hskip 0.1\hsize \hbox to 0.8\hsize{
    \hbox to 22pt{\hfil}#1 \dotfill \hbox to 13pt{\hfil #2}}\par}

  \bigskip
  \centerline {CONTENTS}
  \bigskip
  \baselineskip 13pt

  \ll 1:Introduction:1;
  \ll 2:Approximately proper equivalence relations:2;
  \ll 3:A faithful representation:4;
  \ll 4:Stationary equivalence relations:7;
  \ll 5:Gauge automorphisms:8;
  \ll 6:Finite index:10;
  \ll 7:KMS states:12;
  \ll 8:Existence of KMS states:14;
  \ll 9:Thermodynamic formalism and uniquenes of KMS states:15;
  \ll 10:Conditional minima:20;
  \ll 11:Ground states:22;
  \ll 12:Ground states and maximizing measures:24;
  \ll 13:Phase transitions:27;
  \ni References:28;

  \endgroup

  \section{Introduction}
  An equivalence relation on a compact Hausdorff space is said to be
\stress{proper} when the quotient space is Hausdorff, and
\stress{approximately proper} when it is the union of an increasing
sequence of proper relations.  The first major goal of this paper is
to extend these concepts to \stress{non-commutative spaces}, that is
to C*-algebras, and to construct the corresponding \stress{quotient
space}.  This turns out to be another C*-algebra which is often
non-commutative even when the original algebra is commutative.  An
example of this situation is the \stress{tail-equivalence} relation on
Bernoulli's space whose ``quotient space'' turns out to be the CAR
algebra.

We then introduce the notion of \stress{potentials} and their
associated \stress{gauge actions} which are one-parameter groups of
automorphisms of the ``quotient space''.  A characterization of KMS
states is then provided and we use it to show that KMS states exist
for all values of the inverse temperature.

Starting with a local homeomorphism $T$ on a compact metric space $X$
we consider the equivalence relation on $X$ under which two points $x$
and $y$ are equivalent if there is a natural number $n$ such that
$T^n(x)=T^n(y)$.  This turns out to be an approximately proper
equivalence relation and we apply the abstract theory developed in the
previous sections, enhanced by the use of Ruelle's Perron-Frobenius
Theorem, in order to show uniqueness of KMS states at every
temperature.

  Ground states are studied next and a general characterization of
those states which factor through a certain conditional expectation is
obtained in terms of the support of the corresponding measure.

In the last two sections of the paper we show a relationship between
maximizing measures (the analogue of the Aubry-Mather measures for
expanding maps) and ground states, and  we explore an interesting
example related to phase transitions.

Our construction of the C*-algebra for an approximately proper
equivalence relation should be viewed as a non-commutative
generalization of the groupoid C*-algebra \cite{\RenaultThesis} for
the groupoids treated by Renault in \cite{\RenaultAF} and
\cite{\RenaultRadon}.  In the special case of approximately proper
equivalence relations over commutative algebras, under the assumption
that certain conditional expectations are of index-finite type, an
assumption which we make from section \fcite{6}{} onwards, our
situation actually becomes identical to some situations discussed by
Renault in the above mentioned articles.  Unlike Renault we do not
treat these situations employing groupoids techniques but there is
nevertheless a significant overlap in our conclusions.

  The first named author wishes to acknowledge fruitful discussions
with Jean Renault, Chris Skau, and Anatoly Veshik on topics related to
equivalence relations in the commutative setting.

  \section{Approximately  proper equivalence relations}
  In order to motivate the construction to be made here consider a
compact Hausdorff space $X$ equipped with an equivalence relation $R$.

When the quotient $X/R$ is a Hausdorff space
  we will say that $R$ is a \stress{proper equivalence relation} in
which case the C*-algebra of continuous complex functions on $X/R$,
which we denote as $C(X/R)$, is canonically *-isomorphic to the
subalgebra $C(X;R)$ of $C(X)$ formed by the functions which are
constant on each equivalence class.

On the other hand, given any closed unital *-subalgebra $A\subseteq
C(X)$ define the equivalence relation $R_A$ on $X$ by
  $$
  (x,y)\in R_A \quad \Leftrightarrow\quad \forall f\in A,\ f(x)=f(y).
  $$
  It is then easy to see that $R_A$ is proper and that $C(X;R_A)=A$.
In other words, the correspondence $R \mapsto C(X;R)$ is a bijection
between the set of all proper equivalence relations on $X$ and the
collection of all closed unital *-subalgebras of $C(X)$.

This could be used to give a definition of ``proper equivalence
relations'' over a ``non-commutative space'', that is, a
non-commutative C*-algebra: such a relation would simply be defined to
be a closed unital *-subalgebra.

  This scenario is undoubtedly very neat but it ignores some of the
most interesting equivalence relations in Mathematics, most of which
are not proper.  Consider, for example, the tail-equivalence
relation on Bernouli's space.  The fact that the equivalence
classes are dense implies that $C(X;R)$ consists solely of the
constant functions.  So in this case the subalgebra $C(X;R)$ says
nothing about the equivalence relation we started with.

Fortunately some badly behaved equivalence relations, such as the example
mentioned above, may be described as limits of proper
relations, in the following sense:

  \definition
  An equivalence relation $R$ on a compact Hausdorff space $X$ is said
to be \stress{approximately proper} if there exists an increasing
sequence of proper equivalence relations $\{R_n\}_{n\in\N}$ such that
$R = \bigcup_{n\in\N} R_n$.

We should perhaps say that
  we adopt the convention according to which $\N=\{0,1,2,\ldots\}$.
  Also,
  we view equivalence relations in the strict mathematical sense,
namely as subsets of $X\times X$, hence the set theoretical union
above.

Given $\{R_n\}_{n\in\N}$ as above consider for each $n$ the subalgebra
$\R_n = C(X;R_n)$.  Since $R_n\subseteq R_{n+1}$ we have that
  $\R_n \supseteq \R_{n+1}$.
  Since each $R_n$ may be recovered from $\R_n$ we conclude the
decreasing sequence $\{\R_n\}_{n\in\N}$ encodes all of the information
present in the given sequence of equivalence relations.  We may then
generalize to a non-commutative setting as follows:

  \definition 
  An \stress{approximately proper equivalence relation} on a unital
C*-algebra $A$ is a decreasing sequence $\{\R_n\}_{n\in\N}$ of closed unital
*-subalgebras.  For convenience we will always assume that $\R_0=A$.

It is our goal in this section to introduce a C*-algebra which is
supposed to be the non-commutative analog of the quotient space by an
approximately proper equivalence relation.  A special feature of
our construction is that the resulting algebra is often
non-commutative even when the initial algebra $A$ is commutative.

In order to carry on with our construction it seems that we are
required to choose a sequence of faithful conditional expectations
$\{E_n\}_{n\in\N}$ defined on $A$ with $E_n(A)=\R_n$ and $E_{n+1}\circ E_n =
E_{n+1}$ for every $n$.

Throughout this section, and most of this work, we will therefore fix
a C*-algebra $A$, an approximately proper equivalence relation
$\R=\{\R_n\}_{n\in\N}$, and a sequence $\E=\{E_n\}_{n\in\N}$ of conditional expectations as
above.

  \definition
  \label \DefToep
  The Toeplitz algebra of the pair $(\R,\E)$, denoted $\Toep$, is the
universal C*-algebra generated by $A$ and a sequence
$\{\et_n\}_{n\in\N}$ of projections (self-adjoint idempotents) subject
to the relations:
  \izitem
  \zitem
  $\et_0=1$,
  \zitem $\et_{n+1}\et_n = \et_{n+1}$,
  \zitem
  $\et_n a \et_n = E_n(a) \et_n,$
  \medskip \noindent
  for all $a\in A$ and $n\in\N$.
  
When an element $a\in A$ is viewed in $\Toep$ we will denote it by
$\inc a$. At first glance it is conceivable that the relations above
imply that $\inc a=0$ for some nonzero element $a\in A$.  We will soon
show that this never happens so that we may identify $A$ with its copy
within $\Toep$, and then we will be allowed to drop the underlining notation.

Notice that \lcite{\DefToep.ii} says that the $\et_n$ form a
decreasing sequence of projections.  Also, by taking adjoints in
\lcite{\DefToep.iii}, we conclude that
  $\et_n \inc a \et_n = \et_n \inc{E_n(a)}$ as well.  It follows that each
$\et_n$ lies in the commutant of $\inc {\R_n}$.

  \state Proposition
  \label \AlgFormula
  Given $n,m\in\N$ and $a,b,c,d\in A$ one has that
  $$
  (\inc a\et_n\inc b)(\inc c\et_m\inc d) =
   \cases{
     \inclow {a E_n(b c)}\et_m\inclow d, & if $n\leq m$, \cr
   \cr
     \inclow a \et_n\inclow{E_m(b c)d}, & if $n\geq m$. \cr
   }
  $$

  \proof
  If $n\leq m$ we have
  $$
  (\inc a\et_n\inc b)(\inc c\et_m\inc d) =
  \inc a(\et_n\inc {b c}\et_n)\et_m\inc d =
  \inclow {a E_n(b c)}\et_n\et_m\inclow d =
  \inclow {a E_n(b c)}\et_m\inclow d.
  $$
  If $n\geq m$ the conclusion follows by taking adjoints.
  \proofend 

  \definition
  For each $n\in\N$ we will denote by $\Kt_n$ the closed linear span
of the set $\{\inc a\et_n\inc b: a,b\in A\}$.

  By \lcite{\AlgFormula} we see that for $i\leq n$ one has that both
$\Kt_i \Kt_n$ and $\Kt_n \Kt_i$ are contained in $\Kt_n$.  In
particular each $\Kt_n$ is a C*-subalgebra of $\Toep$.

We now need a  concept borrowed from \scite{\endo}{3.6} and \scite{\tower}{6.2}:

  \definition
  \label \DefineNRedundancy
  Let $n\in\N$.  A finite sequence
  $(k_0,\ldots,k_n)\in \prod_{i=0}^n\Kt_i$
  such that $\sum_{i=0}^n k_ix=0$ for all $x\in \Kt_n$ will be called
an \stress{$n$--redundancy}.  The closed two-sided ideal of $\Toep$
generated by the elements $k_0 + \cdots + k_n$, for all
$n$--redundancies $(k_0,\ldots,k_n)$, will be called the
\stress{redundancy ideal}.

We now arrive at our main new concept:

  \definition
  \label \DefAlg
  The C*-algebra of the pair $(\R,\E)$, denoted $\Alg$, is defined to
be the quotient of $\Toep$ by the redundancy ideal.  Moreover we 
will adopt the following notation:
  \izitem
  \zitem
  The quotient map from $\Toep$ to $\Alg$ will be denoted by $q$.
  \zitem
  The image of $\et_n$ in $\Alg$ will be denoted by $e_n$, 
  \zitem The
image of $\Kt_n$ in $\Alg$ will be denoted by $\K_n$. 

 It is clear that $\K_n$
is the closed linear span of $q(\inc A)e_n q(\inc A)$.

  \section{A faithful representation}
  In this section we will provide a faithful representation of $\Alg$
which will, among other things, show that the natural maps $A\to
\Toep$ and $A\to\Alg$ are injective.

For $n\in\N$ consider the right Hilbert $\R_n$-module $M_n$ obtained by
completing $A$ under the $\R_n$-valued inner product
  $$
  \<a,b\> = E_n(a^*b)
  \for a,b\in A.
  $$
  The canonical map assigning each $a\in A$ to its class in $M_n$ will
be denoted by
  $$
  i_n: A \to M_n.$$
  It is obviously a right $\R_n$-module map.
  For each $a$ in $A$ one may prove that the
  the correspondence
  $$
  i_n(x) \mapsto i_n(ax)
  \for x\in A
  $$
  extends to a map $L^n_a \in \L(M_n)$ (adjointable linear operators
on $M_n$).
  In turn, the correspondence $a\to L^n_a$ may be shown to be an
injective *-homomorphism from $A$ to $\L(M_n)$
  (recall that the $E_n$ are supposed faithful)
  and whenever convenient
we will use it to think of $A$ as subalgebra of $\L(M_n)$.

  We will denote by $\eb_n$ the projection in $\L(M_n)$ obtained by
continuously extending the correspondence $i_n(x) \mapsto i_n(E_n(x))$
to the whole of $M_n$.

Given any two vectors $\xi,\eta\in M_n$ we will denote by
$\Omega_{\xi,\eta}$ the ``generalized rank-one compact operator'' on
$M_n$ given by
  $$
  \Omega_{\xi,\eta}(\zeta) = \xi\<\eta,\zeta\>
  \for \zeta \in M_n.
  $$

  \state Proposition
  \label \EIsOmega
  Given $a,b\in A$ one has that $a \eb_n {b^*} =
\Omega_{i_n(a),i_n(b)}$.
Therefore the closed linear span of the set
  $\{a \eb_n {b^*} : a,b\in A\}$
  is precisely the algebra of generalized compact operators on $M_n$.
This algebra will be denoted  by $\Kb_n$.

  \proof
  For $x\in A$ notice that
  $$
  a \eb_n {b^*}(i_n(x))=
  i_n(a E_n(b^*x)) =
  i_n(a)E_n(b^*x) =
  i_n(a) \<i_n(b),i_n(x)\> =
  \Omega_{i_n(a),i_n(b)}(i_n(x)).
  \proofend
  $$

The following is an important algebraic relation:

  \state Proposition
  \label \CommutRelInRepresentation
  For every $n\in\N$ and every $a\in A$ one has that
  $$
  \eb_n a\eb_n = {E_n(a)}\eb_n = \eb_n{E_n(a)}.
  $$

  \proof
  Given $x\in A$ notice that
  $$
  \eb_n a\eb_n(i_n(x)) =
  i_n\big(E_n(a E_n(x))\big) =
  i_n\big(E_n(a)E_n(x)\big) =
  {E_n(a)}\eb_n(i_n(x)).
  $$
  So
  $\eb_n a\eb_n = {E_n(a)}\eb_n$.  That $\eb_n a\eb_n = \eb_n{E_n(a)}$
follows by taking adjoints.
  \proofend

We now wish to see how do the $M_n$'s relate to each other.

  \state Proposition
  For every $n\in\N$ there exists a continuous $\R_{n+1}$-linear map
  $
  j_n : M_n \to M_{n+1}
  $
  such that
  $j_n(i_n(a)) = i_{n+1}(a)$ for all $a\in A$.
  Moreover for any $\xi,\eta\in M_n$ one has that
  $$
  E_{n+1}(\<\xi,\eta\>) = \<j_n(\xi),j_n(\eta)\>.
  $$

  \proof
  For every $a\in A$ we claim that $\|i_{n+1}(a)\|\leq \|i_n(a)\|$.
In fact
  $$
  \|i_{n+1}(a)\| ^2 =
  \|E_{n+1}(a^*a)\| =
  \|E_{n+1}E_n(a^*a)\| \leq
  \|E_n(a^*a)\| =
  \|i_n(a)\|^2.
  $$
  Thus the correspondence $i_n(a) \mapsto i_{n+1}(a)$ is contractive
and hence extends to a continuous map $j_n:M_n \to M_{n+1}$ such that
$j_n(i_n(a)) = i_{n+1}(a)$.  It is elementary to verify that $j_n$ is
$\R_{n+1}$-linear.  Suppose that $\xi=i_n(a)$ and $\eta =i_n(b)$ where
$a,b\in A$.  Then
  $$
  E_{n+1}(\<\xi,\eta\>) =
  E_{n+1}(\<i_n(a),i_n(b)\>) =
  E_{n+1}(E_n(a^*b)) =
  E_{n+1}(a^*b) =
  \<i_{n+1}(a),i_{n+1}(b)\> \$=
  \<j_{n}(i_n(a)),j_{n}(i_n(b))\> =
  \<j_{n}(\xi),j_{n}(\eta)\>.
  $$
  The conclusion now follows because $i_n(A)$ is dense in $M_n$.
  \proofend
  
  The preceding result gives a canonical relationship between
elements in $M_n$ and $M_{n+1}$.  We will now see how to relate
operators.

  \state Proposition
  There exists an injective *-homomorphism
  $$
  \Phi_n: \L(M_n) \to \L(M_{n+1})
  $$
  such that for $T\in\L(M_n)$ one has that
  $$
  \Phi_n(T)(j_n(\xi)) = j_n(T(\xi))
  \for \xi\in M_n.
  $$

  \proof
  Let $T\in \L(M_n)$.  Since $T^*T \leq \|T\|^2$ one has for all
$\xi\in M_n$ that
  $$
  \<T(\xi),T(\xi)\> =
  \<T^*T(\xi),\xi\> \leq
  \|T\|^2\<\xi,\xi\>.
  $$
  Applying $E_{n+1}$ to the above inequality yields
  $$
  E_{n+1}(\<T(\xi),T(\xi)\>)\leq
  \|T\|^2E_{n+1}(\<\xi,\xi\>),
  $$
  or
  $$
  \<j_n(T(\xi)),j_n(T(\xi))\> \leq
  \|T\|^2(\<j_n(\xi),j_n(\xi)\>),
  $$
  which implies that
  $\|j_n(T(\xi))\| \leq \|T\|\,\|j_n(\xi)\|$. So  the
correspondence
  $$
  j_n(\xi) \mapsto j_n(T(\xi))
  $$
  extends to a bounded linear map $\Phi_n(T):M_{n+1}\to M_{n+1}$ such
that $\Phi(T)(j_n(\xi)) = j_n(T(\xi))$ for all $\xi\in M_n$.

We claim that $\Phi(T)^* = \Phi(T^*)$ for all $T\in \L(M_n)$.  In
order to prove this let $\xi,\eta\in M_n$.  We have that
  $$
  \<j_n(\xi),\Phi(T)(j_n(\eta))\> =
  \<j_n(\xi),j_n(T(\eta))\> =
  E_{n+1}( \<\xi,T(\eta)\> ) \$=
  E_{n+1}( \<T^*(\xi),\eta\> ) =
  \<\Phi(T^*)(j_n(\xi)),j_n(\eta)\>,
  $$
  proving the claim.  It is now easy to see that $\Phi_n$ is indeed a
*-homomorphism from $\L(M_n)$ to $\L(M_{n+1})$.
 
If $T$ is such that $\Phi_n(T)=0$ then for every $\xi\in M_n$ one has
that
  $$
  0 =
  \<\Phi_n(T)(j_n(\xi)),\Phi_n(T)(j_n(\xi))\> =
  \<j_n(T(\xi)),j_n(T(\xi))\> =
  E_{n+1}(\<T(\xi),T(\xi)\>).
  $$
  Since $E_{n+1}$ is faithful we have that $T(\xi)=0$.  Since $\xi$ is
arbitrary we have that $T=0$.
  \proofend

  \definition
  We will denote by $\L_\infty$ the inductive limit of the sequence
  $$
  \L(M_1) \labelarrow{\Phi_1}
  \L(M_2) \labelarrow{\Phi_2}
  \cdots
  $$

Recall that $A$ is viewed as a subalgebra of $\L(M_n)$ via the
correspondence $a \mapsto L^n_a$.  For $a,x\in A$ notice that
  $$
  \Phi_n(L^n_a)(i_{n+1}(x)) =
  \Phi_n(L^n_a)(j_n(i_n(x)) =
  j_n(L^n_a(i_n(x))) =
  j_n(i_n(ax)) =
  i_{n+1}(ax) =
  L^{n+1}_a(i_{n+1}(x)),
  $$
  so that
  $\Phi_n(L^n_a) = L^{n+1}_a$.  It follows that if we identify
$\L(M_n)$ with its image in $\L(M_{n+1})$ under $\Phi_n$ the two
corresponding copies of $A$ will be identified with each other via the
identity map.  Therefore $A$ sits inside of $\L_\infty$ in a canonical
fashion.

  We now claim that
  $\eb_{n+1}\leq \Phi_n(\eb_n)$ for all $n\in\N$.  In fact, for all
$a\in A$
  $$
  \eb_{n+1} \Phi_n(\eb_n) (i_{n+1}(a)) =
  \eb_{n+1} \Phi_n(\eb_n) (j_n(i_n(a))) =
  \eb_{n+1} (j_n(\eb_n(i_n(a)))) =
  \eb_{n+1} (j_n(i_n(E_n(a)))) \$=
  \eb_{n+1} (i_{n+1}(E_n(a))) =
  i_{n+1}(E_{n+1}E_n(a))=
  i_{n+1}(E_{n+1}(a))=
  \eb_{n+1}(i_{n+1}(a)).
  $$

  Within $\L_\infty$ we then get a decreasing sequence of projections
consisting of the images of the $\eb_n$ in the inductive limit, which
we will still denote by $\eb_n$.

  We are now ready to prove the main result of this section whose main
purpose is to give a concrete realization of the so far abstractly
defined  $\Alg$.

  \state Theorem
  \label \ConcreteRepFaithful
  \izitem
  \zitem
  There exists a unique *-homomorphism $\hat\pi:\Toep\to\L_\infty$
such that $\hat\pi(\inc a)=a$ for all $a$ in $A$ and $\hat\pi(\et_n)=\eb_n$
for all $n\in\N$.
  \zitem $\hat\pi$ vanishes on the redundancy ideal and so factors
through $\Alg$ providing a *-homomorphism 
  $$\pi:\Alg\to\L_\infty$$
  such that
  $\pi(e_n)=\eb_n$ and
  $\pi(q(\inc a))=a$,
  where $q$ is the quotient map from $\Toep$ to $\Alg$.
  \zitem $\pi$ is injective and hence $\Alg$ is isomorphic to the
sub-C*-algebra of $\L_\infty$ generated by $A$ and all of the $\eb_n$.

  \proof
  The first point follows from \lcite{\CommutRelInRepresentation}, the
fact that the $\eb_n$ are decreasing, and the universal property of
$\Toep$.

Addressing (ii) all we must show is that $\hat\pi$ vanishes on any element
of the form 
  $$
  s = \sum_{i=0}^n k_i,
  $$
  where $(k_0,\ldots,k_n)$ is an $n$--redundancy. 
  Observing that for $i\leq n$ one has that $\hat\pi(k_i)\in\L(M_i)$ and
that $\L(M_i)$ is contained in $\L(M_n)$ (as subalgebras of the direct
limit $\L_\infty$), we see that $\hat\pi(s)\in\L(M_n)$.  Given $a\in A$
choose $b,c\in A$ such that $E(b^*c)=1$ (e.g.~$b=c=1$) so that
  $$
  \hat\pi(s)\calcat{i_n(a)} = 
  \hat\pi(s)\ \Omega_{i_n(a),i_n(b)}\calcat{i_n(c)} \={\EIsOmega}
  \hat\pi(s)\ (a\eb_n b^*)\calcat{i_n(c)} =
  \hat\pi(s\ \inc a\et_n \inc b^*)\calcat{i_n(c)} = 0
  $$
  because $\inc a\et_n \inc b$ lies in $\Kt_n$.  This shows that
$\hat\pi(s)=0$ and hence proves (ii).

In order to proceed we must now prove that the restriction of $\hat\pi$ to
each $\Kt_n$ is injective.   For this purpose recall from 
\scite{\Watatani}{2.2.9} that $\Kb_n$ is precisely
the \stress{unreduced C*-basic construction} relative to
$E_n$ and thus possesses the universal property described in
\scite{\Watatani}{2.2.7}.
The correspondence 
  $$
  a\in A \mapsto \inc a\in \Toep
  $$
  together with the idempotent $\et_n$ gives by \lcite{\DefToep.iii} a
\stress{covariant representation} of the conditional expectation
$E_n$, according to Definition 2.2.6 in \cite{\Watatani}.  Therefore 
there exists a *-homomorphism $\rho: \Kb_n \to \Kt_n$ such that
$\rho(a\eb_n b)=\inc a\et_n\inc b$ \ for all $a,b\in A$.  It follows that
the composition $\rho \circ\big(\hat\pi|_{\Kt_n}\big)$ is the identity
map hence proving our claim that $\hat\pi|_{\Kt_n}$ is injective.

In order to prove (iii) it suffices to show that for each $n$, $\pi$
is injective on the sub-C*-algebra of $\Alg$ given by
  $$
  B_n = \K_0 + \cdots + \K_n,
  $$
  where the $\K_n$ are defined in \lcite{\DefAlg.iii}
  (note that $B_n$ is indeed a sub-C*-algebra by
\scite{\Ped}{1.5.8}).  In fact, once this is granted we see that $\pi$ 
is isometric on the union of all $B_n$ which is dense in
$\Alg$.  This would prove that $\pi$ is isometric on all of $\Alg$.

Let $b=k_0+\ldots+k_n\in B_n$, where $k_i\in\K_i$, and suppose that
$\pi(b)=0$.  Since $q(\Kt_i)=\K_i$ we may write $k_i=q(\hat k_i)$,
where the $\hat k_i\in \Kt_i$.  We therefore have that
  $\hat\pi(\hat k_0+\cdots+\hat k_n)=0$.

  We now claim that $(\hat k_0,\ldots,\hat k_n)$ is an
$n$-redundancy.  In order to prove it let $x\in\Kt_n$ and note that 
  $(\hat k_0+\cdots+\hat k_n) x \in \Kt_n$ by  \lcite{\AlgFormula}.
But since 
  $\hat\pi\big((\hat k_0+\cdots+\hat k_n) x\big)=0$ and $\hat\pi$ is
injective on $\Kt_n$ we have that $(\hat k_0+\cdots+\hat k_n) x=0$ as
claimed.  So $\hat k_0+\cdots+\hat k_n$ lies in the redundancy ideal
and hence
  $$
  b = 
  q(\hat k_0+\cdots+\hat k_n)=0.
  \proofend
  $$

  \state Corollary
  The maps 
  $$
  a\in A\to \inc a\in \Toep
  $$
  and
  $$a\in A\to q(\inc a)\in
\Alg
  $$
  are injective.

  \proof 
  Follows immediately from our last result.
  \proofend

From now on we will therefore identify $A$ with  $\inc A$ and also with
$q(\inc A)$.  

  \section{Stationary equivalence relations}
  \label\StationarySect
  In this section we will study approximately proper equivalence
relations which have a specially simple description.

  \definition An approximately proper equivalence relation
$\R=\{\R_n\}_{n\in\N}$ over a unital C*-algebra $A$ is said to be
\stress{stationary} if there exists an unital injective *-endomorphism
$\a: A \to A$ such that $\R_{n+1}=\a(\R_n)$ for all $n$.

In this case observe that $\R_n$ is simply the range of $\a^n$.  Throughout
this section we will fix a stationary approximately proper equivalence
relation $\R=\{\R_n\}_{n\in\N}$ over $A$.  We will also fix an
endomorphism $\a$ as above.

Let $E$ be a given faithful conditional expectation from $A$ to
$\R_1$.  Define conditional expectations $E_n$ from $A$ to $\R_n$
by
  $$
  E_n = \a^{n-1} \underbrace{(E \a\inv) \ldots (E
\a\inv)}_{n-1\rm\;times}E.
  $$
  It is easy to see that $E_{n+1}\circ E_n = E_{n+1}$ for every $n$.

  \definition We will say that a sequence of conditional
expectations $\E=\{E_n\}_{n\in\N}$ is stationary if it is obtained as above
from a single faithful conditional expectation $E:A\to\R_1$.

  Throughout this section we will fix a stationary sequence of
conditional expectations as above.  Observe that the composition 
$\Tr=\a\inv E$ is a \stress{transfer operator} in the sense of
\scite{\endo}{2.1}.  We may then form the crossed-product $\CP$ as in
\scite{\endo}{3.7}.  Denote by $\gamma$ the \stress{scalar gauge
action} \scite{\tower}{3.3} on $\CP$.

  The main result we wish to present in this section is in order:

  \state Theorem
  \label\TwoConstructions
  With the hypothesis introduced in this section  $\Alg$ is isomorphic to
the sub-C*-algebra of the crossed-product algebra $\CP$ formed by the
fixed points for the scalar gauge action.

  \proof
  Follows immediately from \scite{\tower}{6.5} since the algebra
$\check{\cal U}$ mentioned there (see also \scite{\tower}{4.8}) is
isomorphic to $\Alg$ by \lcite{\ConcreteRepFaithful.iii}.  For the
proof that the fixed point algebra is precisely $\check{\cal U}$ see
\scite{\M}{4.1}.
  \proofend

We may now finally give a nontrivial example of our construction.  Let
$A$ be an $n\times n$ matrix of zeros and ones without any zero rows or
columns and let $(X,T)$ be the corresponding Markov sub-shift.  Define
the endomorphism $\a$ of $C(X)$ by $\a(f) = f\circ T$, for all $f$ in
$C(X)$.  Also let $E$ be the conditional expectation from $C(X)$ to
the range of $\a$ given by
  $$
  E(f)\calcat x  =
  {1 \over  \#\{y: T(y)=x\}} \displaystyle \sum_{T(y)=x} f(y)
  \for f \in C(X),\ x\in X.
  $$

  We may then form $\R$ and $\E$ as above.  By \scite{\endo}{6.2} one
has that $C(X)\mathop{\crossproduct_{\a,\Tr}} {\bf N}$ is the
Cuntz-Krieger algebra ${\cal O}_A$.  By \lcite{\TwoConstructions} we
then have that $\Alg$ is isomorphic to the subalgebra of ${\cal O}_A$
formed by the fixed point algebra for the gauge action.  When $n=2$
and $A = \pmatrix{1 & 1\cr 1 & 1}$ we then have that $\Alg$ is
isomorphic to the CAR algebra.

  \section{Gauge automorphisms}
  In this section we return to the general case, therefore fixing  a
C*-algebra $A$, an approximately proper equivalence relation
$\R=\{\R_n\}_{n\in\N}$, and a sequence $\E=\{E_n\}_{n\in\N}$ of
compatible conditional expectations as before.

  We wish to  introduce the notions of \stress{potentials}
and their corresponding \stress{gauge automorphisms} which will be the object
of study of later sections.  We start with a simple technical fact:

  \state Proposition
  \label \ENKcommute
  Let $i\leq n$ and let $b\in\Ctr(\R_i)$ (meaning the center of\/ $\R_i$) then
  $$
  E_n(a b) = E_n(b a)
  \for a\in A.
  $$
  
  \proof
  We have
  $$
  E_n(a b) =
  E_n\big(E_i(a b)\big) =
  E_n\big(E_i(a)b\big) =
  E_n\big(b E_i(a)\big) =
  E_n\big(E_i(b a)\big) =
  E_n(b a).
  \proofend
  $$

  \definition
  \label \DefinePotential
  By a \stress{potential} we will mean a sequence
  $z = \{z_n\}_{n\in\N}$
  such that  $z_n$ belongs to $\Ctr(\R_n)$ for every $n\in\N$.

Given a potential $z$ observe that every $z_n$ commutes with every
other $z_m$.  Therefore we may set
  $$
  z^{[n]} = z_0 z_1 \ldots z_{n-1}
  \for n\in\N,
  $$
  without worrying about the order of the factors.  We will also use
the notation $z^{-[n]}$ to mean $\(z^{[n]}\)\inv$ when the latter exists.

If $w = \{w_n\}_{n\in\N}$ is another potential it is clear that 
$z w := \{z_n w_n\}_{n\in\N}$ is again a potential and that
  $$
  (z w)^{[n]}= z^{[n]} w^{[n]}
  \for n\in\N.
  $$

  Potentials may be used to define automorphisms as follows:

  \state Proposition
  Let $u = \{u_n\}_{n\in\N}$ be a unitary potential (in the sense that
each $u_n$ is a unitary element).  Then there is an automorphism
$\tphi_u$ of\/ $\Toep$ such that
  $$
  \tphi_u(a) = a
  \for a\in A,
  $$
  and
  $$
  \tphi_u(\et_n) = u^{[n]}\et_n u^{-[n]}
  \for n\in\N.
  $$
  Moreover,  given another unitary potential $v$ one has that
$\tphi_{u v} =
\tphi_u\tphi_v$.
  
  \proof
  For every $n\in\N$ let $f_n=u^{[n]}\et_n u^{-[n]}$.  Then
  $$
  f_n f_{n+1} =
  u^{[n]}\et_n u^{-[n]} u^{[n+1]}\et_{n+1} u^{-[n+1]} =
  u^{[n]}\et_n u_n\et_{n+1} u^{-[n+1]} \$= 
  u^{[n]}u_n\et_n \et_{n+1} u^{-[n+1]} =
  u^{[n+1]}\et_{n+1} u^{-[n+1]} =
  f_{n+1},
  $$
  so the $f_n$ are decreasing.  For $a\in
A$ we have
  $$
  f_n a f_n =
  u^{[n]}\et_n u^{-[n]} a u^{[n]}\et_n u^{-[n]} =
  u^{[n]}E_n\(u^{-[n]} a u^{[n]}\)\et_n u^{-[n]} =
  $$$$
  \=\ENKcommute
  u^{[n]}E_n(a)\et_n u^{-[n]} =
  E_n(a)u^{[n]}\et_n u^{-[n]} =
  E_n(a)f_n.  
  $$
  By the universal property of $\Toep$ there exist a *-homomorphism
$\tphi_u:\Toep\to\Toep$ satisfying the conditions in the statement,
except possibly for the fact that $\tphi_u$ is an automorphism.

Given another unitary potential $v$ one can easily prove that
$\tphi_{u v} = \tphi_u\tphi_v$ by checking on the generators.
Plugging $v=u\inv:=\{u_n\inv\}_{n\in\N}$ we than have that
$\tphi_{u\inv}$ and $\tphi_u$ are each others inverse and hence
$\tphi_u$ is an automorphism.
  \proofend

  \state Proposition
  For every unitary potential $u$ the automorphism $\tphi_u$ leaves
the redundancy ideal invariant (in the sense that 
the image of
the redundancy ideal under $\tphi_u$ is exactly the redundancy
ideal)
and hence drops to an automorphism
$\varphi_u$ of\/ $\Alg$ which is the identity on $A$ and such that 
  $$
  \varphi_u(e_n) = u^{[n]}e_n u^{-[n]}
  \for n\in\N.
  $$

  \proof
  It is elementary to verify that $\tphi_u(\Kt_n)=\Kt_n$ for all $n$.
  Thus, if 
  $(k_0,\ldots,k_n)$  
  is a  redundancy we have that 
  $(\tphi_u(k_0),\ldots,\tphi_u(k_n)) \in \prod_{i=0}^n\Kt_i$.
Moreover if $x\in\Kt_n$ we have that
  $$
  \sum_{i=0}^n \tphi_u(k_i)x =
  \tphi_u\(\sum_{i=0}^n k_i \varphi_u\inv(x)\) = 0.
  $$
  Therefore
  $(\tphi_u(k_0),\ldots,\tphi_u(k_n))$ is a redundancy and hence 
  $\tphi_u(k_0+\cdots+k_n)$ lies in the redundancy ideal.  So we see
that $\tphi_u$ sends the redundancy ideal in itself.  Since the same
holds for $\tphi_{u\inv} = \tphi_u\inv$ if follows that the image of
the redundancy ideal under $\tphi_u$ is precisely the redundancy
ideal and hence the proof is concluded.
  \proofend

So far we have introduced single gauge automorphisms, but now we would like
to define one-parameter groups of such:

  \definition
  \label\DefGauge
  \izitem
  \zitem
  A potential $h = \{h_n\}_{n\in\N}$ is said to be strictly positive
when for each $n$ there exists a real number $c_n>0$ such that
$h_n\geq c_n$.
  \zitem Given a strictly positive potential $h = \{h_n\}_{n\in\N}$
and a complex number $z$ we denote by $h^{z}$ the potential
$\{h_n^z\}_{n\in\N}$, and by $h^{z[n]}=(h^z)^{[n]}$, for $n\in\N$.
  \zitem
  The \stress{gauge action} for a strictly positive potential $h$ is
the one-parameter group $\sigma = \{\sigma_t\}_{t\in\Reals}$ of
automorphisms of $\Alg$ given by
  $\sigma_t = \varphi_{h^{it}}$
  for all $t\in\Reals$.

  \medskip
  Given $a,b\in A$ and $n\in\N$  observe that
  $$
  \sigma_t(a e_n b) = a h^{it[n]} e_n h^{-it[n]} b
  \for a,b\in A
  \for n\in\N.
  \eqno{(\seqnumbering)}
  \label \SigmaOnMonomial
  $$
  It is therefore clear that the gauge action is strongly continuous.

  \section{Finite index}
  Starting with this section we will restrict ourselves to the case in
which the $E_n$ are of index-finite type according to
  \scite{\Watatani}{1.2.2}.
  We refer the reader to \cite{\Watatani} for
the basic definitions and facts about index-finite type conditional
expectations, which will now acquire a preponderant role in our study. 

  \state Proposition
  \label \RestrictionFiniteIndex
  If $E_m$ is of index-finite type then its restriction to each $\R_n$,
where $n\leq m$, is also of index-finite type.  Moreover if\/
$\{u_1,\ldots,u_k\}$ is a quasi-basis for $E_m$ then
$\{E_n(u_1),\ldots,E_n(u_k)\}$ is a quasi-basis for the restriction of
$E_m$ to $\R_n$.

  \proof
  For every $a\in \R_n$ we have that
  $$
  a = 
  E_n(a) =
  E_n\Big(\sum_{i=0}^k u_i E_m(u_i^*a)\Big) =
  \sum_{i=0}^k E_n(u_i) E_m\big(E_n(u_i^*a)\big) =
  \sum_{i=0}^k E_n(u_i) E_m\big(E_n(u_i)^*a\big).
  \proofend
  $$

  \state Proposition
  \label \IncreasingKn
  Let $n\leq m$.  Suppose that the restriction of $E_m$ to $\R_n$ is
of index-finite type and let $\{v_1,\ldots,v_k\}\subseteq \R_n$ be a quasi-basis for
it.  Then
  \izitem
  \zitem $\sum_{i=0}^k v_i e_m v_i^* = e_n$,
  \zitem $\K_n\subseteq \K_m$.

  \proof
  Let $a,b\in A$ and observe that
  $$
  \Big(\et_n - \sum_{i=0}^k v_i \et_m v_i^*\Big)a \et_m b =
  \et_n a \et_n  \et_m b  - \sum_{i=0}^k v_i E_m(v_i^*a) \et_m b =
  E_n(a) \et_m b - \sum_{i=0}^k v_i E_m(E_n(v_i^*a)) \et_m b \$=
  E_n(a) \et_m b - \sum_{i=0}^k v_i E_m(v_i^*E_n(a)) \et_m b =  
  E_n(a) \et_m b - E_n(a) \et_m b = 0.
  $$
  Therefore the ($m+1$)-tuple
  $$
  \Big(0,\ldots,0, \et_n,0,\ldots,0,-\sum_{i=0}^k v_i \et_m v_i^*\Big)
  $$
  is an $m$-redundancy from where (i) follows.  Obviously (ii) follows from
(i).
  \proofend

  \state Corollary
  If all of the $E_n$ are of index-finite type then the $\K_n$ are
increasing and $\Alg$ is the closure of\/ $\bigcup_{n\in\N} \K_n$.

  \proof
  By \lcite{\RestrictionFiniteIndex} we have that $E_{n+1}|_{\R_n}$ is
of index-finite type.  Hence by \lcite{\IncreasingKn} we have that
$\K_n\subseteq\K_{n+1}$.  Since $A = \K_0$ and  for every $n$
we have that
  $
  e_n \in \K_n
  $
  the conclusion follows.
  \proofend

In the finite index case we have the following elementary description
of the $\K_n$:

  \state Proposition
  If all of the $E_n$ are of index-finite type then $M_n=i_n(A)$ and $\K_n =
L_{\R_n}(A)$, where $L_{\R_n}(A)$ denotes the set of all (not
necessarily adjointable or even continuous) additive right
$\R_n$-linear maps on $A$ (where $A$ is identified with $M_n$ via
$i_n$).

  \proof By \scite{\Watatani}{2.1.5} there exists a constant
$\lambda_n>0$ such that
  $\|E_n(a^*a)\|^{1/2} \geq \lambda_n \|a\|$, for all $a$ in $A$.
  Therefore
  $$
  \|i_n(a)\| =
  \|E_n(a^*a)\|^{1/2} \geq
  \lambda_n \|a\|,
  $$
  so that $i_n$ is a Banach space isomorphism onto its range which is
therefore a complete normed space, hence closed.  Since $i_n(A)$ is
dense in $M_n$ we conclude that $i_n(A)=M_n$.  We will therefore
identify $M_n$ and $A$.

It is clear that $\K_n \subseteq L_{\R_n}(A)$.  In order to prove the
converse inclusion let
$\{u_1,\ldots,u_n\}$ be a quasi-basis for $E_n$.
  Then, given any additive $\R_n$-linear map $T$ on $A$ and $a\in A$ we
have
  $$
  T(a) =
  T\Big(\sum_{i=1}^m u_i E_n(u_i^* a)\Big) =
  \sum_{i=1}^m T(u_i) E_n(u_i^* a) =
  \sum_{i=1}^m T(u_i) \<u_i,a\> =
  \sum_{i=1}^m \Omega_{T(u_i),u_i}(a),
  $$
  so that
  $
  T = \sum_{i=1}^m \Omega_{T(u_i),u_i} \in \Kb_n.
  $
  \proofend

  This last result gives a curious description of the dense subalgebra
$\bigcup_{n\in\N} \K_n$ of $\Alg$, namely that it is formed by the additive
operators which are linear with respect to \stress{some} $\R_n$.
Observe that this is not quite the same as requiring linearity with
respect to the intersection of the $\R_n$!

  One of the main tools in our study from now on will be a certain
conditional expectation from $\Alg$ to $A$.  Unfortunately we can only
show its existence in the finite-index case.

  \state Proposition
  \label \MainCondExpectation
  If all of the $E_n$ are of index-finite type then there exists a
conditional expectation 
  $$
  G : \Alg \to A,
  $$
  such that for each $n\in\N$ one has that
  $$
  G(e_n) = \lambda_0\inv \ldots \lambda_{n-1}\inv,
  $$
  where $\lambda_n=\ind (E_{n+1}|_{\R_n})$.
  If $A$ is commutative then $G$ is the unique conditional expectation 
from $\Alg$ to $A$.

  \proof
  Set $\lambda_n=\ind (E_{n+1}|_{\R_n})$ so that
  $\lambda = \{\lambda_n\}_{n\in\N}$ is a potential in the sense of
definition \lcite{\DefinePotential} and the proposed value for
$G(e_n)$ above is just $\lambda^{-[n]}$.  Observe moreover that
$\lambda^{-[n]}$ commutes with $\R_{n-1}$.

Let $n\in\N$ be fixed.  Observing that $\K_n$ is isomorphic to $\Kb_n$
by \lcite{\ConcreteRepFaithful.iii} and arguing exactly as in
\scite{\tower}{8.4} we conclude that there exists a positive
$A$--bimodule map
  $G_n : \K_n \to A$ such that
  $G_n(e_n) = \lambda^{-[n]}$.

  We claim that $G_{n+1}$ extends $G_n$.  In fact  let
$\{u_1,\ldots,u_k\}$ be a quasi-basis for $E_{n+1}$.  Then by
\lcite{\RestrictionFiniteIndex} we have that
$\{E_n(u_1),\ldots,E_n(u_k)\}$ is a quasi-basis for $E_{n+1}|_{\R_n}$.

  By \lcite{\IncreasingKn.i} we have that
  $
  e_n =   \sum_{i=1}^k E_n(u_i) e_{n+1} E_n(u_i)^*,
  $
  so that
  $$
  G_{n+1}(e_n) =
  \sum_{i=1}^k E_n(u_i) \lambda^{-[n+1]} E_n(u_i)^* =
  \lambda^{-[n+1]}  \sum_{i=1}^k E_n(u_i) E_n(u_i)^* \$=
  \lambda^{-[n+1]}  \ind (E_{n+1}|_{\R_n}) =
  \lambda^{-[n+1]}  \lambda_n =
  \lambda^{-[n]} =
  G_n(e_n).
  $$
  The claim then follows easily from the fact that both $G_n$ and
$G_{n+1}$ are $A$--bimodule maps. 

  As a consequence we see that each $G_n$ restricts to the identity on
$A$ and hence $G_n$ is a conditional expectation from $\K_n$ to $A$.
Conditional expectations are always contractive so  there exists a
common extension $G:\Alg\to A$ which is the desired map.
  
Suppose that $A$ is commutative and that $G'$ is another conditional
expectation from $\Alg$ to $A$.  Given $n$ let $\{u_1,\ldots,u_k\}$ be
a quasi-basis for $E_n$ and hence by \lcite{\IncreasingKn.i} we have
  $$
  1 = G'(1) =
  G'\Big(\sum_{i=0}^k u_i e_n u_i^*\Big) =
  \sum_{i=0}^k u_i G'(e_n) u_i^* =
  G'(e_n) \ind(E_n),
  $$
  so necessarily $G'(e_n) = \ind(E_n)\inv = \lambda^{-[n]}$ by
\scite{\Watatani}{1.7.1}).  Once knowing that $G$ and $G'$ coincide on
the $e_n$ it is easy to see that $G= G'$.
  \proofend

  \section{KMS states}
  In this section we will begin the general study of KMS states for
gauge actions on $\Alg$.  We refer the reader to \cite{\BR} and
\cite{\Ped} for the basic theory of KMS states.

Given what are probably limitations in our methods we will all but
have to assume that $A$ is commutative.  To be precise we will suppose
from now on that the conditional expectations $E_n$ satisfy the
following trace-like property:
  $$
  E_n(a b)=E_n(b a)\for a,b\in A,
  \eqno{(\seqnumbering)} 
  \label\TraceLikeProperty
  $$
  which is obviously the case when $A$ is commutative.  Unfortunately
we have no interesting non-commutative example of this situation but
since we do not really have to suppose that $A$ is commutative and in
the hope that some such example will be found we will proceed without
the commutativity of $A$.

We will moreover assume that all of the $E_n$ are of index-finite type
and will denote by $G$ the conditional expectation given by
\lcite{\MainCondExpectation}.  Our first result is that any KMS state
factors through $G$.

  \state Proposition
  \label\KMSFactors
  Let $h$ be a strictly positive potential, let $\b>0$, and let $\phi$ be a
$(\sigma,\b)$-KMS state (i.e.~a KMS state for $\sigma$ at inverse
temperature $\b$) on $\Alg$ for the gauge action $\sigma$ associated
to $h$.  Then $\phi = \phi\circ G$.

  \proof
  Given $a,b\in A$ and $n\in\N$ it is clear from
\lcite{\SigmaOnMonomial} that $a e_n b$ is an analytic element with
  $$
  \sigma_z(a e_n b) = a h^{i z[n]} e_n h^{-i z[n]} b
  \for z\in {\bf C}.
  $$
  We claim that
  $$
  \phi(a e_n b) = \phi\(h^{\b[n]} E_n(b a h^{-\b[n]})e_n\)
  \for a,b\in A\for n\in \N.
  \eqno{(\dagger)}
  $$
  In order to prove it we use the KMS condition
  as follows
  $$
  \phi(a e_n b) = 
  \phi(e_n b a) = 
  \phi(e_n b a \sigma_{i\b}(e_n)) = 
  \phi(e_n b a h^{-\b[n]}e_n h^{\b[n]}) \$= 
  \phi(E_n(b a h^{-\b[n]})e_n h^{\b[n]}) = 
  \phi(h^{\b[n]} E_n(b a h^{-\b[n]})e_n),
  $$
  proving $(\dagger)$.  We next claim that 
  $$
  \phi(a e_{n+1}) = \phi(\lambda_n\inv a e_n) 
  \for a\in A,
  $$
  where $\lambda_n$ is defined in \lcite{\MainCondExpectation}.  In
order to prove this claim let
  $\{v_1,\ldots,v_k\}\subseteq \R_n$ be a quasi-basis for the
restriction of $E_{n+1}$ to $\R_n$.  Then by \lcite{\IncreasingKn.i}
we have for all $x\in A$ that
  $$
  \phi(x e_n) = 
  \phi\Big( \sum_{i=0}^k x v_i e_{n+1} v_i^*\Big) \= {\dagger}
  \sum_{i=0}^k  \phi\(h^{\b[n+1]}E_{n+1}(v_i^*x v_i
h^{-\b[n+1]})e_{n+1}\).
  $$
  Since $v_i\in \R_n$ and since $h^{-\b[n+1]}$ commutes with $\R_n$ we
have that
  $$
  E_{n+1}(v_i^*x v_i h^{-\b[n+1]}) =
  E_{n+1}(v_i^*x h^{-\b[n+1]} v_i ) =
  E_{n+1}(x h^{-\b[n+1]} v_i v_i^*),
  $$
  by the trace-like property of $E_{n+1}$.
  We then conclude that
  $$
  \phi(x e_n) = 
  \phi\(h^{\b[n+1]} E_{n+1}(x h^{-\b[n+1]} \lambda_n)e_{n+1}\).
  $$
  Using $(\dagger)$ once more we have that
  $$
  \phi(a e_{n+1}) = \phi\(h^{\b[n+1]} E_{n+1}(a h^{-\b[n+1]})e_{n+1}\).
  $$  
  So when $x=\lambda_n\inv a$ we have that
  $\phi(x e_n) = \phi(a e_{n+1})$ which is precisely the identity
we were looking for.
  By induction we then have that
  $$
  \phi(a e_n) =
  \phi(\lambda^{-[n]}a).
  $$
  Therefore for all $a,b\in A$
  $$
  \phi(a e_n b) =
  \phi(b a e_n ) =
  \phi(\lambda^{-[n]}b a) =
  \phi(a \lambda^{-[n]}b ) =
  \phi(G(a e_n b)).
  $$  
  As the closed linear span of the set of elements of the form
$a e_n b$ is dense in $\Alg$ the proof is complete.
  \proofend

  In particular it follows that every KMS state is determined by its
restriction to $A$.  It is therefore useful to know which states on
$A$ occur as the restriction of a KMS state.

  \state Proposition
  \label \MyKMSCond
  Let $\phi$ be a state on $A$ and let $\b>0$.  Then the composition
$\psi = \phi\circ G$ is a $(\sigma,\b)$-KMS state if and only if
  $$
  \phi(a) =
  \phi\(\Lambda^{-[n]} E_n(\Lambda^{[n]} a)\)
  \for a\in A
  \for n\in\N,
  $$
  where $\Lambda=\{\Lambda_n\}_{n\in\N}$ is the potential given by
$\Lambda_n = h_n^{-\b} \lambda_n$.

  \proof 
  Suppose that $\psi$ is a $(\sigma,\b)$-KMS state.  Then for all
$a,b,c,d\in A$ and all $n\in \N$ we have
  $$
  \psi\((a e_n b)\sigma_{i\b}(c e_n d)\) =
  \psi\((c e_n d)(a e_n b)\).
  \eqno{(\dagger)}
  $$  
  Observe that the left hand side of $(\dagger)$ equals
  $$
  \psi\(a e_n b c h^{-\b[n]} e_n h^{\b[n]} d\) =
  \psi\(a E_n(b c h^{-\b[n]}) e_n h^{\b[n]} d\) =
  \phi\(a E_n(b c h^{-\b[n]}) \lambda^{-[n]}  h^{\b[n]} d\).
  $$
  Meanwhile the right hand side of $(\dagger)$ equals
  $$
  \psi\(c E_n(d a) e_n b\) =
  \phi\(c E_n(d a) \lambda^{-[n]} b\).
  $$
  Plugging $b=1$, $c = h^{\b[n]}$, and $d = h^{-\b[n]} \lambda^{[n]}$
we have that $(\dagger)$ implies that
  $$
  \phi(a) =
  \phi\( h^{\b[n]} E_n( h^{-\b[n]} \lambda^{[n]} a) \lambda^{-[n]}\) =
  \phi\(\Lambda^{-[n]} E_n(\Lambda^{[n]} a)\).
  $$

  In order to prove the converse we first claim that if $\phi$
satisfies the condition in the statement for $n=1$ then $\phi$ must be
a trace.  
In fact,  observing that
  $
  \Lambda^{[1]} = \Lambda_0 \in \Ctr(A)
  $
  we have for all $a,b\in A$ that
  $$
  \phi(a b) =
  \phi\(\Lambda^{-[1]} E_1(\Lambda^{[1]} a b)\) =
  \phi\(\Lambda^{-[1]} E_1( a\Lambda^{[1]} b)\) =
  \phi\(\Lambda^{-[1]} E_1(  \Lambda^{[1]} b a)\) =
  \phi(b a),
  $$
  where we have again used the trace-like property of $E_1$.
  Supposing now that $\phi$ satisfies the above condition not only for
$n=1$ but for all $n\in\N$ let us prove that $\psi$ is a KMS state.  For
this we would like to prove that
  $$
  \psi\((a e_n b)\sigma_{i\b}(c e_m d)\) =
  \psi\((c e_m d)(a e_n b)\),
  \eqno{(\ddagger)}
  $$
  for all $a,b,c,d\in A$ and $n,m\in \N$.  Supposing that $n\leq m$
the left hand side of $(\ddagger)$ equals
  $$
  \psi\(a e_n b c h^{-\b[m]} e_m h^{\b[m]} d\) =
  \psi\(a E_n(b c h^{-\b[m]}) e_m h^{\b[m]} d\) \$=
  \phi\(a E_n(b c h^{-\b[m]}) \lambda^{-[m]} h^{\b[m]} d\) =
  \phi\(E_n(b c h^{-\b[m]}) h^{\b[m]} \lambda^{-[m]} d a\) = \ldots
  $$
  Letting $x = h_n^\b\ldots h_{m-1}^\b$ observe that $x\in \R_n$ and 
$h^{\b[m]} = x h^{\b[n]}$ so the above equals
  $$
  \ldots =
  \phi\(E_n(b c h^{-\b[n]}x\inv) x h^{\b[n]} \lambda^{-[m]} d a\) =
  \phi\(E_n(b c h^{-\b[n]}) h^{\b[n]} \lambda^{-[m]} d a\) \$=
  \phi\(\Lambda^{-[n]} E_n\(\Lambda^{[n]} E_n(b c h^{-\b[n]})
h^{\b[n]} \lambda^{-[m]} d a \) \) =
  \phi\(\Lambda^{-[n]} E_n(b c h^{-\b[n]}) E_n(\lambda^{[n]}
\lambda^{-[m]} d a ) \).
  $$
  Meanwhile the right hand side of $(\ddagger)$ equals
  $$
  \psi\(c e_m E_n(d a) b\) =
  \phi\(c \lambda^{-[m]} E_n(d a) b\) =
  \phi\(b c \lambda^{-[m]} E_n(d a) \) \$=
  \phi\(\Lambda^{-[n]} E_n\(\Lambda^{[n]} b c \lambda^{-[m]} E_n(d a)
\)\) =
  \phi\(\Lambda^{-[n]} E_n\(b c \lambda^{-[m]} \lambda^{[n]} h^{-\b[n]} \) E_n(d a)\).
  $$
  Observing that $\lambda^{-[m]} \lambda^{[n]} \in \R_n$ we therefore
see that $(\ddagger)$ is proved under the hypothesis that $n\leq m$.
If, on the other hand, $n\geq m$ the left hand side of $(\ddagger)$
becomes
  $$
  \psi\(a e_n b c h^{-\b[m]} e_m h^{\b[m]} d\) =
  \psi\(a e_n E_m(b c h^{-\b[m]}) h^{\b[m]} d\) \$=
  \phi\(a \lambda^{-[n]}  E_m(b c h^{-\b[m]}) h^{\b[m]} d\) =
  \phi\(
    \Lambda^{-[m]} E_m\(\Lambda^{[m]} 
    \lambda^{-[n]}  E_m(b c h^{-\b[m]}) h^{\b[m]} d a
    \)\) \$=
  \phi\(\Lambda^{-[m]} E_m(b c h^{-\b[m]}) E_m\( h^{\b[m]} d a
\Lambda^{[m]} \lambda^{-[n]} \)\) =
  \phi\(\Lambda^{-[m]} E_m(b c h^{-\b[m]}) E_m\( d a
\lambda^{[m]} \lambda^{-[n]} \)\).
  $$
  The right hand side of $(\ddagger)$ equals
  $$
  \psi\(c E_m(d a)e_n b\) =
  \phi\(c E_m(d a)\lambda^{-[n]} b\) =
  \phi\(\lambda^{-[n]} b c E_m(d a) \) \$=
  \phi\( \Lambda^{-[m]} E_m\(\Lambda^{[m]} \lambda^{-[n]} b c E_m(d
a)\) \)=
  \phi\( \Lambda^{-[m]} E_m\(b c h^{-\b[m]} \lambda^{[m]} \lambda^{-[n]} \) E_m(d
a) \).
  $$
  The conclusion follows once more because 
$\lambda^{[m]} \lambda^{-[n]} \in \R_m$.
  \proofend

Putting together our last two results we reach one of our main
goals:

  \state Theorem
  \label\MainCharacterization
  Let $\R$ be an approximately proper equivalence relation on a
C*-algebra $A$ and let $\E=\{E_n\}_{n\in\N}$ be a sequence of conditional
expectations of index-finite type defined on $A$ with $E_n(A)=\R_n$
satisfying \lcite{\TraceLikeProperty} and $E_{n+1}\circ E_n = E_{n+1}$
for every $n$.  Also let $h$ be any strictly positive potential and
denote by $\sigma$ the associated gauge action on $\Alg$.  Then for
every $\b>0$ the correspondence $\psi\mapsto\phi=\psi|_A$ is a bijection
from the set of $(\sigma,\b)$-KMS states $\psi$ on $\Alg$
and the set of states $\phi$ on $A$ satisfying
  $$
  \phi(a) =
  \phi\(\Lambda^{-[n]} E_n(\Lambda^{[n]} a)\)
  \for a\in A
  \for n\in\N,
  $$
  where $\Lambda=\{\Lambda_n\}_{n\in\N}$ is the potential given by
$\Lambda_n = h_n^{-\b} \lambda_n$.  The inverse of this correspondence
is given by $\phi\mapsto\psi=\phi\circ G$, where $G$ is given in
\lcite{\MainCondExpectation}.

  \section{Existence of KMS states}
  Theorem \lcite{\MainCharacterization} gives a precise
characterization of the KMS states on $\Alg$ in terms of states on
$A$ satisfying certain conditions.  It does not say, however, if such
states exist.  We will now take up the task of showing the existence
of at least one KMS state for each inverse temperature $\b>0$.  We
begin with a technical result which says that the conditions on $\phi$
required by \lcite{\MyKMSCond} increase in strength with $n$.

  \state Proposition
  \label \NPlusOneImpliesN
  Let $\phi$ be a state on $A$ and suppose that the formula
  $$
  \phi(a) =
  \phi\(\Lambda^{-[n]} E_n(\Lambda^{[n]} a)\)
  \for a\in A,
  $$
  holds for $n=k+1$, where $k\in\N$ is given.  Then the formula holds for
$n=k$.

  \proof
  For each $n\in\N$ let $F_n$ be the operator on $A$ given by 
  $$
  F_n(a) = \Lambda^{-[n]} E_n(\Lambda^{[n]} a)
  \for a\in A.
  $$
  Then the formula in
the statement is equivalent to $F_n^*(\phi)=\phi$, where $F_n^*$
refers to the transpose operator on the dual of $A$.

We claim that for all $n$ one has that
  $F_{n+1}\circ F_n = F_{n+1}$.  In fact,
  observing that $\Lambda^{[n+1]} 
  \Lambda^{-[n]} = \Lambda_n\in \R_n$, we have
  $$ 
  F_{n+1}\(F_n(a)\) =
  \Lambda^{-[n+1]} E_{n+1}\(\Lambda^{[n+1]} 
  \Lambda^{-[n]} E_n(\Lambda^{[n]} a)  \) \$=
  \Lambda^{-[n+1]} E_{n+1}\(E_n( \Lambda^{[n+1]} 
  \Lambda^{-[n]} \Lambda^{[n]} a)  \) =
  \Lambda^{-[n+1]} E_{n+1}( \Lambda^{[n+1]}a)) =
  F_{n+1}(a).
  $$
  Given that $F_{k+1}^*(\phi)=\phi$ we have
  $$
  F_k^*(\phi) =
  F_k^*(F_{k+1}^*(\phi)) =
  (F_{k+1}F_k)^*(\phi) =
  F_{k+1}^*(\phi) =
  \phi.
  \proofend
  $$

We now arrive at the main result of this section.

  \state Theorem
  \label \ExistenceKMS
  Let $\R$ be an approximately proper equivalence relation on a
C*-algebra $A$ and let $\E=\{E_n\}_{n\in\N}$ be a sequence of conditional
expectations of index-finite type defined on $A$ with $E_n(A)=\R_n$
satisfying \lcite{\TraceLikeProperty} and $E_{n+1}\circ E_n = E_{n+1}$
for every $n$.  Also let $h$ be any strictly positive potential and
denote by $\sigma$ the associated gauge action on $\Alg$.  Then for
every $\b>0$ there exists at least one $(\sigma,\b)$-KMS state on $\Alg$.
  
  \proof
  \def\S{{\curly S}}
  For each $n\in\N$ let $\S_n$ be set of all states on $A$ satisfying
$F_n^*(\phi)=\phi$, where $F_n$ is the operator defined in the
beginning of the proof of \lcite{\NPlusOneImpliesN}.  It is clear that
the $\S_n$ are closed subsets of the state space of $A$ and hence
compact.

  We claim that $\S_n$ is nonempty for every $n$.  In order to prove
this let $\tau$ be any trace on $A$.  Observe that traces on $A$ may
be obtained by composing any state with $E_1$.  For a given $n$, let
$\phi = F_n^*(\tau)$.  Since $F_n^2=F_n$ it is clear that
$F_n^*(\phi)=\phi$.  Moreover $\phi$ is a positive linear functional
because for all $a\in A_+$ we have
  $$
  \phi(a) =
  \tau\(\Lambda^{-[n]} E_n(\Lambda^{[n]} a)\) =
  \tau\(\Lambda^{-{1\over 2}[n]} E_n\(\Lambda^{{1\over 2}[n]} a
\Lambda^{{1\over 2}[n]}\) \Lambda^{-{1\over 2}[n]}\) \geq 0.
  $$
  Dividing $\phi$ by $\phi(1)$ 
  (observe that $\phi(1)\neq 0$ by \scite{\Watatani}{2.1.5}) 
  thus gives an element of $\S_n$ so that
$\S_n\neq\emptyset.$ By \lcite{\NPlusOneImpliesN} we have that the
$\S_n$ are decreasing so their intersection is nonempty.  Any $\phi$
belonging to that intersection is a state on $A$ satisfying the
condition in \lcite{\MyKMSCond} and hence $\phi\circ G$ is a
$(\sigma,\b)$-KMS state on $\Alg$.
  \proofend

It should be noticed that the method employed above may be used to
give an iterative process to produce KMS states: start with any
state $\phi_0$ on $A$ and define
  $$
  \phi_{n}= \phi_{n-1}(F_n(1))\inv F_n^*(\phi_{n-1}).
  $$
  Any weak accumulation point of the sequence $\{\phi_n\}_n$ will be a
state $\phi$ on $A$ satisfying \lcite{\MyKMSCond} and hence $\phi\circ
G$ is the desired KMS state.

In the present level of generality there is not much more we can say
about KMS states.  In the following sections we will discuss an
example in which KMS states will be proven to be unique as well.

  \section{Thermodynamic formalism and uniquenes of KMS states}
  \label \ThermoSection
  In this part of the paper we will show a relationship between the
KMS states we have been discussing and the Gibbs states of
Thermodynamic Formalism, as developed by Bowen, Ruelle, and
Sinai \cite{\Bo}, \cite{\Rone}, \cite{\Rdois}, \cite{\Rtres}.

Throughout the rest of this section we will fix a compact metric space
$X$ and a local homeomorphism $T:X \to X$.  We will also let $\a$ be
the endomorphism of $C(X)$ given by
  $$
  \a(f) = f\circ T
  \for f\in C(X).
  $$

  Consider the
equivalence relation on $X$ given by
  $$
  x \sim y \quad \Leftrightarrow \quad \exists n\in \N,\ T^n(x) =T^n(y).
  $$
  In the case of the left shift on Bernouli's space (an example to be
kept in the back of one's mind) this equivalence relation turns out to
be the tail-equivalence relation which is not proper.  However it is
easy to see that it is always  approximately proper, and that it is the union
of the equivalence relations $R_n$ given by
  $$
  (x,y)\in R_n \quad \Leftrightarrow \quad T^n(x) =T^n(y).
  \eqno{(\seqnumbering)} 
  \label\RnForTransformation
  $$
  Clearly each $R_n$ is proper and the algebra $C(X;R_n)$ is precisely
the range of $\a^n$.  For simplicity we will denote the latter algebra
by $\R_n$.

  We now need conditional expectations $E_n$ from $C(X)$ onto $\R_n$
and these will be obtained as follows.  By the assumption that $T$ is
a local homeomorphism and that $X$ is compact we see that $T$ is
necessarily a covering map.  The inverse image under $T$ of each $x\in
X$ is therefore a finite set.  Given a continuous strictly positive
function $p:X\to \Reals$ consider the associated
Ruelle-Perron-Frobenius operator given by
  $$
  \Tr_p (f)\calcat x = \sum_{T(z)=x} p(z) f(z)
  \for f\in C(X),\ x\in X.
  $$
  We will assume that $p$ is such that $\Tr_p$ is normalized (meaning
that $\Tr_p(1)=1$). This means that for every $x\in X$ the
association $z \mapsto p(z)$ is a probability distribution on the
equivalence class of $x$ relative to $R_1$.

  It is easy to show that $\Tr_p$ satisfies the identity
  $$
  \Tr_p(f)g = \Tr_p(f\a(g))
  \for f,g\in C(X).
  \eqno{(\seqnumbering)}
  \label\MyTransferCondition
  $$
  For any $n\in\N$ set
  $$
  E_n =\a^n\Tr_p^n.
  \eqno{(\seqnumbering)}
  \label\EsFromL
  $$
  Given $f\in C(X)$ one then has that $E_1(f)\calcat
x$ is just the weighted average of $f$ over the equivalence class of $x$
relative to $R_1$.  Therefore $E_1$ is a conditional expectation onto
$\R_1$.  Likewise $E_n$ is a conditional expectation onto $\R_n$ and
because the composition $\Tr_p\circ\a$ is the identity map on $C(X)$
we have that
  $E_m\circ E_n = E_m$ for $m \geq n$.  Setting  $\R=\{\R_n\}_{n\in\N}$
and  $\E=\{E_n\}_{n\in\N}$ we may then speak of $\Alg$.

Observe that the present situation is precisely that of a stationary
equivalence relation described in section \lcite{\StationarySect}.

Given any $f\in C(X)$ it is clear that $\a^n(f)\in\R_n$ for all $n$
and hence the sequence $\{\a^n(f)\}_{n\in\N}$ is a potential.
Accordingly we will adopt the notation $f^{[n]}$ to mean
  $$
  f^{[n]} = f \a(f)\ldots \a^{n-1}(f).
  $$

For later use it is convenient to give an explicit description for
$\Tr_p^n$ as well as $E_n$:

  \state Lemma
  \label \ExplicitEn
  Let $n\in\N$ then for every $f\in C(X)$ and $x\in X$ one has that
  $$
  \Tr_p^n(f)\calcat x =
  \sum_{T^n(z)=x} p^{[n]}(z) f(z),
  $$
  and 
  $$
  E_n(f)\calcat x =
  \sum_{(z,x)\in R_n} p^{[n]}(z) f(z),
  $$

  Before giving the proof we should notice that in summations of the form
$\displaystyle \sum_{(z,x)\in R_n}$, which will be often used
from now on, the variable which we mean to sum  upon will always be the first
one mentioned ($z$ in this case) even though equivalence
relations are well known to be symmetric.

  \proof (of \ExplicitEn)
  In order to prove the first statement we use induction on $n$
observing that the case $n=1$ follows by definition.   Given $n\geq 1$ we
have
  $$
  \Tr_p^{(n+1)}(f)\calcat x =
  \Tr_p^n(\Tr_p(f))\calcat x =
  \sum_{T^n(z)=x} p^{[n]}(z) \sum_{T(w)=z} p(w) f(w) = \ldots
  $$
  Notice that a pair $(z,w)$ is such that $T^n(z)=x$ and $T(w)=z$ if
and only if it is of the form $(T(w),w)$ where
$T^{n+1}(w)=x$.  Therefore the above equals
  $$
  \ldots =
  \sum_{T^{n+1}(w)=x} p^{[n]}(T(w)) p(w) f(w)  =
  \sum_{T^{n+1}(w)=x} p^{[n+1]}(w) f(w),
  $$
  proving the first statement.  The second statement then follows
easily.
  \proofend

In the sequel we  compute the index of our conditional expectations.

  \state Proposition
  \label \OurCondOfFiniteTyp
  For each $n\in N$ we have that $E_{n+1}|_{\R_n}$ is of index-finite
type and $\ind(E_{n+1}|_{\R_n}) = \a^n(p\inv)$.

  \proof
  Let $\{V_i\}_{i=1}^m$ be a finite open covering of $X$ such that the
restriction of\/ $T$ to each $V_i$ is one-to-one and let
$\{v_i\}_{i=1}^m$ be a partition of unit subordinate to this covering.
Set $u_i=(p\inv v_i)^{1/2}$ and observe that for every $f\in
C(X)$ and $x\in X$ one has that
  $$
  \sum_{i=1}^m u_i E_1(u_if)\calcat x =
  \sum_{i=1}^m u_i(x) 
    \sum_{\buildrel {\scriptstyle z\in X} \over {T(z)=T(x)}}
p(z) u_i(z)f(z) \$=
  \sum_{i=1}^m u_i(x) p(x) u_i(x)f(x)=
  \sum_{i=1}^m v_i(x) f(x) = f(x).
  $$
  Therefore $\{u_1,\ldots,u_m\}$ is a quasi-basis for $E_1$ so that
  $$
  \ind(E_1) =
  \sum_{i=1}^m u_i^2 =
  \sum_{i=1}^m p\inv v_i = p\inv.
  $$
  Next  observe that the  diagram
  $$
  \matrix{ \R_0 & \map{E_1} & R_1\cr \cr
           {\scriptstyle \a^n} \downarrow \quad &&
           \qquad \downarrow {\scriptstyle \a^n} \cr\cr
           \R_n & \map{E_{n+1}} & \R_{n+1}}
  $$
  is commutative.  Therefore $E_{n+1}|_{\R_n}$ is conjugate to $E_1$
under $\a^n$ and so $\ind (E_{n+1}|_{\R_n}) = \a^n(\ind(E_1)) =
\a^n(p\inv)$.
  \proofend

We therefore have that each $E_n$ is of index-finite type.  Also
notice that in the notation of \lcite{\MainCondExpectation} we have
proven that $\lambda_n= \a^n(p\inv)$.

Let $H$ be a strictly positive continuous function on $X$.
Setting $h_n = \a^n(H)$ for every $n\in\N$ we have that
$h:=\{h_n\}_{n\in\N}$ is a strictly positive potential in the sense of
\lcite{\DefGauge}.  The corresponding gauge action will be denoted by
$\sigma$.

We are interested in showing that for every $\b>0$ there exists a
unique $(\sigma,\b)$-KMS state on $\Alg$, thus improving on Theorem
\lcite{\ExistenceKMS}.  

 Given $\b>0$ consider the Ruelle-Perron-Frobenius operator
associated to $H(z)^{-\beta}$, namely
  $$
  \Tr\hb (f) \calcat x = \sum_{T(z)=x} H(z)^{-\beta}  f(z)
  \for f\in C(X),\ x\in X.  
  $$

In order to achieve our goal we will need to use the celebrated
Ruelle-Perron-Frobenius Theorem whose conclusions are:

  \sysstate{Conclusions of the Ruelle-Perron-Frobenius}{\sl}{
  \label \ConclusionsRPF
  \item{a)}
There exists a unique pair $(c\hb, \nu\hb)$ such that $c\hb$ is a
strictly positive real number, $\nu\hb$ is a probability measure on
$X$, and
  $$\Tr\hb^{*} (\nu\hb) = 
  c\hb \nu\hb,
  $$
  where $\Tr\hb^{*}$ refers to the transpose operator on the dual of
$C(X)$, which in turn is
identified with the space of finite regular Borel measures on $X$.
  \item{b)} There exists a strictly positive continuous function
$k\hb$ on $X$ such that
  \medskip  \itemitem{$\bullet$}
  $\displaystyle\int k\hb \d \nu\hb  =1$,
  \medskip  \itemitem{$\bullet$}
   $  \Tr\hb(k\hb)= \ev\hb k\hb$, and
  \medskip  \itemitem{$\bullet$}
  $\displaystyle  \lim_{n \to \infty} {\Tr\hb^n(f) \over \ev\hb^n} =
  k\hb \int f \d \nu\hb 
  \for f\in C(X),
  $
  \bigskip
  where the limit is with respect to the (sup) norm topology of
$C(X)$.  
 }

 Initally proven for the shift on the one-sided Bernouli's space
\scite{\Rone}{Theorem 3} this Theorem has been proved to hold under
more general hypothesis: see for example
  \cite{\Bo},          
  \cite{\Rdois},       
  \cite{\Wone},        
  \cite{\Craizer},     
  \cite{\Rtres},       
  \cite{\Fan},         
  \cite{\Keller},      
  \cite{\Baladi},      
  \cite{\FanJiangOne}, 
  \cite{\Wtwo},        
  \cite{\FanJiangTwo}. 

  The reader is referred to the above articles for more details on the
various hypothesis under which the Ruelle-Perron-Frobenius holds so we
will simply assume its conclusions as above.

Later we will consider a situation in which the above conclusions do
not hold causing the phenomena of phase-transitions. This model is
known as the Fisher-Felderhof model \cite{\Ltwo}, \cite{\Lthree},
\cite{\FL}.  

  \definition
  The probability $\nu\hb$ is called the \stress{Gibbs state}
associated to  $H^{-\b}$.

In the sequel we show the following elementary relationship between
the operators $\Tr_p$ and $\Tr\hb$:

  \state Proposition
  Given $\b>0$ and  $n\in\N$ we have that
  $$
  \Tr\hb^n(f) = \Tr_p^n(\Lambda^{[n]}f)
  \for f\in C(X),
  $$
  where the potential 
$\Lambda=\{\Lambda_n\}_{n\in\N}$ was defined in 
\lcite{\MyKMSCond} by $\Lambda_n
= h_n^{-\b} \lambda_n$.
  
  \proof In the present situation we have that $h_n=\a^n(H)$ and
$\lambda_n= \a^n(p\inv)$ so that
  $$
  \Lambda_n =
  \a^n(H)^{-\b} \a^n(p\inv) =
  \a^n(H^{-\b} p\inv).
  $$
  Next observe that for $f\in C(X)$ we have
  $$
  \Tr\hb(f) = \Tr_p(H^{-\b} p\inv  f) = \Tr_p(\Lambda_0 f).
  $$
  The conclusion now follows easily  by
induction using  \lcite{\MyTransferCondition}.
  \proofend

We will now show that the Gibbs states indeed give KMS states on $\Alg$:

  \state Proposition
  \label\GibsIsKMS
  For every $\b>0$ the state $\phi\hb$ on $C(X)$ corresponding via the
Riesz representation Theorem to the Gibbs state $\nu\hb$
satisfies the conditions of \lcite{\MyKMSCond} and hence the
composition $\psi\hb = \phi\hb\circ G$ is a $(\sigma,\b)$-KMS state on $\Alg$.

  \proof
  The condition that $\nu\hb$ is an eigenmeasure for $\Tr\hb$ 
gives for every $f\in C(X)$ and any $n\in\N$ that
  $$
  \phi\hb(\Tr_p^n(\Lambda^{[n]}f)) = 
  \phi\hb(\Tr\hb^n(f)) = 
  c\hb^n \phi\hb(f).
  $$
  Plugging $f=\Lambda^{-[n]}\a^n(g)$  above, where $g\in C(X)$, we obtain
  $$
  \phi\hb(g) = c\hb^n \phi\hb(\Lambda^{-[n]}\a^n(g)).
  $$
  In order to prove the condition in
\lcite{\MyKMSCond} we then compute
  $$
  \phi\hb\(\Lambda^{-[n]} E_n(\Lambda^{[n]} f)\) =
  \phi\hb\(\Lambda^{-[n]} \a^n \Tr_p^n(\Lambda^{[n]} f)\) =
  c\hb^{-n}  \phi\hb\(\Tr_p^n(\Lambda^{[n]} f)\) =
  \phi\hb(f).
  $$
  This concludes the proof.
  \proofend  

Our next main goal will be to show that the state $\psi\hb$ given by
the above result is the unique $(\sigma,\b)$-KMS state on $\Alg$. 

  \state Theorem
  \label\MainUniqueness
  Let $T$ be a local homeomorphism on a compact metric space $X$ and
consider the approximately proper equivalence relation
$\R=\{R_n\}_{n\in\N}$, where each $R_n$ is given by
\lcite{\RnForTransformation}. Let $p:X\to \Reals$ be a strictly
positive continuous function satisfying $\sum_{T(z)=x}p(z)=1$ for
every $x\in X$ and define the sequence of conditional expectations
$\E=\{E_n\}_{n\in\N}$ as in \lcite{\EsFromL}.  Let $H$ be a strictly
positive continuous function on $X$ and consider the one parameter
automorphism group of $\Alg$ given by the potential $h:=\{H\circ
T^n\}_{n\in\N}$.  Assuming \lcite{\ConclusionsRPF} we have that for
every $\b>0$ the state $\psi\hb$ given by \lcite{\GibsIsKMS} is the
unique $(\sigma,\b)$-KMS state on $\Alg$.

  \proof
  Let $\psi$ be a $(\sigma,\b)$-KMS state on $\Alg$ and let $\phi$ be its
restriction to $C(X)$.  By \lcite{\KMSFactors} we have that
$\psi=\phi\circ G$ so it suffices to show that $\phi=\phi\hb$.  Fix
$f\in C(X)$ and notice that by \lcite{\MyKMSCond} we have
  $$
  \phi(f) = 
  \phi\(\Lambda^{-[n]} E_n(\Lambda^{[n]} f)\) =
  \phi\(\Lambda^{-[n]} \a^n\Tr_p^n(\Lambda^{[n]} f)\) =
  \phi\(\Lambda^{-[n]} \a^n\Tr\hb^n(f)\) \$=
  \ev\hb^n  \phi\(\Lambda^{-[n]} \a^n\({\Tr\hb^n(f) \over
  \ev\hb^n}\)\).
  \eqno{(\dagger)}
  $$
  We next
claim that if we replace the argument of $\a^n$ in $(\dagger)$ by its
limit, namely $\phi\hb(f) k\hb$, we will arrive at an expression which
converges to $\phi(f)$ as $n\to\infty$.  In order to prove this we compute
  $$
  \left|\phi(f) - \ev\hb^n \phi\(\Lambda^{-[n]} \a^n\(\phi\hb(f)
k\hb\)\)\right| =
  \left| \ev\hb^n  \phi\(\Lambda^{-[n]} \a^n\({\Tr\hb^n(f) \over
  \ev\hb^n} - \phi\hb(f) k\hb\)\)\right| \$ \leq
  \ev\hb^n  \phi\(\Lambda^{-[n]}\) \[{\Tr\hb^n(f) \over
  \ev\hb^n} - \phi\hb(f) k\hb\].
  $$
  The claim will be proven once we show that the expression $\ev\hb^n
\phi(\Lambda^{-[n]})$ is bounded from above with $n$.  In fact, as $k\hb$ is
strictly  positive, there exists $m>0$ such that $k\hb>m$.
Therefore plugging $f:=k\hb$ in $(\dagger)$ leads to
  $$
  \phi(k\hb) = 
  \ev\hb^n  \phi\(\Lambda^{-[n]} \a^n\(k\hb\)\) \geq
  \ev\hb^n  \phi\(\Lambda^{-[n]} \) m,
  $$
  from where one easily deduces the desired boundedness.  
  Summarizing we have proven that
  $$
  \phi(f) =
  \phi\hb(f) \lim_{n\to\infty} \ev\hb^n\ \phi\(\Lambda^{-[n]}
\a^n\(k\hb\)\),
  $$
  for every $f\in C(X)$.  Since both $\phi$ and $\phi\hb$ evaluate to
$1$ on the constant function $f=1$, it follows that
$\displaystyle\lim_{n\to\infty} \ev\hb^n\ \phi\(\Lambda^{-[n]}
\a^n\(k\hb\)\) = 1$ and hence that $\phi=\phi\hb$ as desired.
  \proofend

  As a consequence we have:

  \state Corollary
  Let $X$, $T$, $\R$, $p$, and $\E$ be as in \lcite{\MainUniqueness}.
Then $\Alg$ admits a unique trace.

  \proof
  Set $H=1$ in \lcite{\MainUniqueness} so that the corresponding one
parameter automorphism group is the trivial one.  Fixing an arbitrary
$\b>0$ observe that the $(\sigma,\b)$-KMS states on $\Alg$ are
precisely the traces.  The conclusion then follows
from \lcite{\MainUniqueness}.
  \proofend

  \section{Conditional minima}
  So far we have studied KMS states at positive temperature and we
have seen how they relate to the Gibbs states of statistical
mechanics.  We next want to discuss ground states but before that we
need to study the notion of conditional minimum points.

Our discussion in this and the next section may be viewed as a special
case of Renault's study of ground-state cocycles over groupoids
\scite{\RenaultThesis}{Section 3}.  We begin with some notation:

  \definition 
  \label \NotationForM
  Let 
  $R$ be a proper equivalence relation on the compact space $X$,
  let $h$ be a  continuous real function on $X$, and 
  let $C$ be a closed subset of $X$.
  We denote by:
  \izitem
  \zitem
  $M_{h,C}$ the set of minimum points for $h$ over $C$, namely
  $$
  M_{h,C} = \left\{ x\in C: h(x) = \inf_{y\in C}h(y)\right\},
  $$
  \zitem
  $M_h^R$ the union of all $M_{h,C}$ as $C$ range in the quotient
space $X/R$ (observe that each $C \in X/R$ is a closed subset of $X$).

  \bigskip
  Observe that a necessary and sufficient condition for $x$ to be in
$M_h^R$ is that
  $$
  \forall y\in X,\ (x,y)\in R \Rightarrow h(x) \leq
h(y).
  \eqno{(\seqnumbering)}
  \label \CondForMemberMn
  $$
  For this reason the points in  $M_h^R$ should  be called
conditional minimum points of $h$.
  Observe also that our hypotheses imply that $M_{h,C}$ is nonempty for
every $C\in X/R$ so one sees that $M_h^R$ meets every single equivalence class.

Even though $M_{h,C}$ is closed for every equivalence class $C$ it may
be that $M_h^R$ is not closed.  However under suitable conditions we
may assure that $M_h^R$ is closed:

  \state Proposition
  {\rm (see \scite{\RenaultThesis}{3.16.iii})}
  \label \MClosed
  Let 
  $R$ be a proper equivalence relation on the compact space $X$ and
  let $h$ be a continuous real function on $X$.  If $R$ is open
(recall that an equivalence relation is said to be open when the
saturation of each open set is open) then $M_h^R$ is closed.

  \proof
  Let $\pi$ be the quotient map from $X$ to the (Hausdorff) space
$X/R$.  Observe that the hypothesis that $R$ is open implies that
$\pi$ is an open mapping.

Let $\{x_i\}_i$ be a net in $M_h^R$ converging to a point $x$ in
$X$ which we assume by contradiction does not belong to $M_h^R$.
Therefore there exists $y$ in $X$ with $(x,y)\in R$ such that
$h(y)<h(x)$.  Let $\alpha$ be any real number with
  $$
  h(y) < \alpha < h(x),
  $$
  and let $U$ be the open set given by
  $$
  U = \{t\in X: h(t) < \alpha\},
  $$
  so that $y\in U$. Observe that
  $$\pi(x_i) \labelarrow{n\to\infty} \pi(x) = \pi(y) \in \pi(U).$$
Since $\pi$ is an open mapping we have that $\pi(U)$ is open so there
exists some $i_0$ such that for all $i\geq i_0$ one has that
$\pi(x_i)\in \pi(U)$.  Given that $h(x)>\alpha$ there exists $j \geq
i_0$ such that $h(x_j)>\alpha$.
  Insisting that $\pi(x_j)\in \pi(U)$ pick $z\in U$ such that
$\pi(z)=\pi(x_j)$ and observe that because $z\in U$ we have that
  $$
  h(z) < \alpha < h(x_j),
  $$
  which contradicts the fact that $x_j\in M_h^R$.
  \proofend

So far we have been considering a proper equivalence relation $R$ on a
compact set $X$ and a continuous real function $h$ on $X$.  From now
on we will assume that $R$ is such that the quotient map $\pi:X\to
X/R$ is a covering map, which incidentally implies that $R$ is open.
We wish to add to this setup a conditional expectation $E$ from $C(X)$
to $\R:=C(X;R)$ which will be obtained as follows: fix a strictly
positive continuous function $p$ on $X$ and let $E:C(X)\to \R$ be given
by
  $$
  E(f)\calcat x = \sum_{(y,x)\in R} p(y) f(y)
  \for f\in C(X),\ x\in X.
  \eqno{(\seqnumbering)}
  \label \OneMoreConExp
  $$
  If we assume  that
  $$
  \sum_{(y,x)\in R} p(y) =1
  \for x\in X
  \eqno{(\seqnumbering)}
  \label \PSumsOne
  $$ 
  it is easy to see that $E$ is indeed a conditional
expectation onto $\R$.

  The following is the main result of this section.  It will be the
crucial technical tool in our characterization of ground states.

  \state Lemma
  \label \TechnicalGround
  Let $R$ be a proper equivalence relation on a compact space $X$ such
that the corresponding quotient map is a covering map.  Let $p$ be a
strictly positive continuous function on $X$ satisfying
\lcite{\PSumsOne} and define the conditional expectation $E$ as in
\lcite{\OneMoreConExp}.  If $h$ is another strictly positive
continuous functions on $X$ define for each real number $\b\geq 0$
the operator $E^\b$ on $C(X)$  by
  $$
  E^\b(f) = h^\b E\(h^{-\b} f\)
  \for f\in C(X).
  $$
  Then for every  probability measure $\mu$ on $X$ the following conditions
are equivalent:
  \izitem
  \zitem The support of $\mu$ is contained in $M_h^R$,
  \zitem For every $f,g\in C(X)$ one has that
  $$
  \sup_{\b\geq 0} \left | \int_X f E^\b(g) \d \mu \right |\leq
\|f\|\, \|g\|,
  $$
  \zitem For every $f,g\in C(X)$ one has that
  $$
  \sup_{\b\geq 0} \left | \int_X f E^\b(g) \d \mu \right |< \infty,
  $$
  \zitem The inequality in (iii) holds for $f = g = 1$.
  
  \proof We begin with the proof that (i) implies (ii).  For this let
$f,g\in C(X)$, and $\b\geq 0$.  We have by  (i)  that
  $$
  \int_X f E^\b(g) \d \mu =
  \int_{M_h^R}f(x) h^\b(x) E(h^{-\b}g)\calcat x \d \mu(x) \$=
  \int_{M_h^R}f(x)
h^\b(x)\sum_{(y,x)\in R}p(y)h^{-\b}(y)g(y) \d \mu(x)
\$=
  \int_{M_h^R}\sum_{(y,x)\in R} p(y)
  \({ h(x) \over h(y)}\)^\b
  f(x)g(y) \d \mu(x).
  $$
  For $x\in M_h^R$ and $y$ such that $(y,x)\in R$ we have by
\lcite{\CondForMemberMn} that   $\({h(x) \over h(y)}\)^\b\leq 1$.
Therefore
  $$
  \left | \int_X f E^\b(g) \d \mu \right | \leq
  \int_{M_h^R}\sum_{(y,x)\in R} p(y) \big|f(x) g(y)\big| \d \mu(x) \leq
  \|f\|\,\|g\|.
  $$

  It is evident that (ii) $\Rightarrow$ (iii) $\Rightarrow$ (iv) so it
remains to prove that (iv) $\Rightarrow$ (i).  For this purpose assume
(iv) and suppose that $x_0\in X\setminus M_h^R$.  It suffices to show
that there exists a neighborhood $U$ of $x_0$ such that $\mu(U)=0$.

Given that $x_0\notin M_h^R$ there exists $y_0\in X$ such that
$(y_0,x_0)\in R$ and $h(x_0) > h(y_0)$.  One may then
choose a real number $c>1$, and open sets $U$ and $V$ with $x_0\in U$,
$y_0\in V$, and such that $h(x) > c h(y)$ whenever $x\in
U$ and $y\in V$.  By reducing the size of both $U$ and $V$ we may
assume that $\pi$, the quotient map, is a homeomorphism restricted to each of $U$ and $V$
and also such that $\pi(U)=\pi(V)$.  Consequently there exists a
homeomorphism $\tau : U \to V$ such that $\pi(x) = \pi(\tau(x))$ for all
$x\in U$.  In particular
  $$
  {h(x) \over h(\tau(x))} > c
  \for x\in U.
  $$
  
Let $\ds K = \sup_{\b\geq 0} \int_X E^\b(1) \d \mu$, which is finite
by (iv), so that for every $\b\geq0$ one has
  $$
  K \geq
  \int_X E^\b(1) \d \mu =
  \int_{X}\sum_{(y,x)\in R} p(y)\({ h(x) \over
h(y)}\)^\b \d \mu(x) \geq
  \int_{U} p(\tau(x))\({ h(x) \over h(\tau(x))}\)^\b \d
\mu(x),
  $$
  where the last inequality is a consequence of replacing $X$ by the
smaller set $U$ and replacing the sum by a single summand, namely when
$y = \tau(x)$. Let $m$ be the (positive) infimum of $p$ over $X$ so
that we conclude that
  $$
  K \geq
  m c^\b \mu(U),
  $$
  or equivalently that
  $\mu(U) \leq c^{-\b}m\inv K$.  Since $\b$ is arbitrary and $c>1$ we have that
$\mu(U)=0$ as desired.
  \proofend
  
  \section{Ground states}
  In this section we will apply the conclusions reached above to study
ground states on $\Alg$.  The setup for now will be as follows: $X$
will be a compact Hausdorff space and $\R=\{R_n\}_{n\in\N}$ an approximately
proper equivalence relation on $X$.  We will also fix a real potential
$h = \{h_n\}_{n\in\N}$.  Recall from \lcite{\DefinePotential} that
this means that each $h_n$ is a continuous real functions
in $\R_n:= C(X;R_n)$.

  \state Proposition
  \label  \MnDecreasing
  For every $n\in\N$ let $M_n$ be the set of conditional minimum
points of $h^{[n]}$ relative to $R_n$, namely 
  $$
  M_n = {M\vrule height 10pt depth 5pt width 0pt}_{h^{[n]}}^{R_n}
  $$
  in the notation of \lcite{\NotationForM.ii}.
  Then $M_{n+1}\subseteq M_n$.

  \proof
  Let $x\in M_{n+1}$.  In order to show that $x\in M_n$ we will employ
the characterization given in \lcite{\CondForMemberMn}.  So let $y$ be
such that $(x,y)\in R_n$.  Since the $R_k$ are increasing we have that
$(x,y)\in R_{n+1}$ and hence 
  $$
  h^{[n+1]}(x) \leq   h^{[n+1]}(y). 
  \eqno{(\dagger)}
  $$
  Observe that because $h_n$ belongs to $C(X;R_n)$ we have that
$h_n(x)=h_n(y)$.  Dividing both sides of $(\dagger)$ by this common
value leads to $h^{[n]}(x) \leq   h^{[n]}(y)$, completing the proof.
  \proofend

If one tries to apply the definition of conditional minimum points for
the relation $R = \bigcup_{n\in\N} R_n$, which we are attempting to
approximate by the sequence $\{R_n\}_{n\in\N}$, one is likely to run
into some trouble, not least because equivalence classes need not
always be closed (in fact they are often dense).  An alternative approach
is to look at points which are conditional minima for all of the
$R_n$.

  \definition 
  Given an approximately proper equivalence relation $\R=\{R_n\}_{n\in\N}$ on
a compact space $X$ and a real potential $h = \{h_n\}_{n\in\N}$ we
will denote by $M_h^\R$ the intersection of the $M_{h^{[n]}}^{R_n}$ as
$n$ range in $\N$.

  Observe that if all of the $R_n$ are open equivalence relations it
follows from \lcite{\MClosed} and \lcite{\MnDecreasing} that $M_h^\R$
is a nonempty compact subset of $X$.  

From this point on we will assume that the $R_n$ are not only open but
also that the quotient maps are covering maps as in
\lcite{\TechnicalGround}.  In addition to this we will fix a strictly
positive potential $p = \{p_n\}_{n\in\N}$.  Following
\lcite{\ExplicitEn} and \lcite{\OneMoreConExp} we define maps
  $E_n:C(X)\to \R_n$ by
  $$
  E_n(f)\calcat x = \sum_{(y,x)\in R_n} p^{[n]}(y) f(y)
  \for f\in C(X),\ x\in X.
  $$
  
  \state 
  \label \ConditionOnPForCondexp
  Lemma
  Suppose that for every $n$ and every $R_{n+1}$-equivalence class $C$ one
has that 
  $$
  \sum_{D} p_n(D) = 1,
  $$
  where the sum extends over all $R_n$-equivalence classes $D$
contained in $C$, and for each such $D$ one interprets $p_n(D)$ as the
common value of $p_n(x)$ for any $x\in D$.  Then each $E_n$ is a
conditional expectation of index-finite type onto $\R_n$ and
$E_{n+1}\circ E_n = E_{n+1}$.

  \proof
  We first claim that for every $n\in\N$ and every $x\in X$ one has that
$\ds\sum_{(y,x)\in R_n} p^{[n]}(y) = 1$.  In order to prove it we use
induction observing that the case ``$n=1$'' follows from the hypothesis.
Assuming that $n\geq1$ we have
  $$
  \sum_{(y,x)\in R_{n+1}} p^{[n+1]}(y) =
  \sum_{i=1}^n \sum_{y\in C_i} p^{[n+1]}(y) = \cdots
  $$
  where $\{C_1,\ldots,C_n\}$ is the decomposition of the
$R_{n+1}$-equivalence class of $x$ into $R_n$-equivalence classes.
  The above then equals
  $$
  \cdots = 
  \sum_{i=1}^n \sum_{y\in C_i} p_n(y)p^{[n]}(y) =
  \sum_{i=1}^n p_n(C_i) \sum_{y\in C_i} p^{[n]}(y) =
  \sum_{i=1}^n p_n(C_i) = 1,
  $$
  where the penultimate equality follows from the induction hypothesis
and the last equality is a consequence of  our hypothesis.
  It immediately follows that $E_n$ is in fact a conditional
expectation onto $\R_n$.  The proof that $E_n$ is of index-finite type
is a simple modification of \lcite{\OurCondOfFiniteTyp} and hence will
be omitted.

With respect to the last  part of the
statement let $f\in C(X)$ so that for $x\in X$ we have
  $$
  E_{n+1} (E_n(f))\calcat x =
  \sum_{(y,x)\in R_{n+1}} p^{[n+1]}(y)  \sum_{(w,y)\in R_n} p^{[n]}(w)
f(w) = \cdots
  $$
  Letting $\{C_1,\ldots,C_n\}$ be as in the first part of the proof we
have that the above equals
  $$
  \cdots =
  \sum_{i=1}^n \sum_{y\in C_i} p^{[n+1]}(y)  \sum_{w\in C_i} p^{[n]}(w)
f(w) =
  \sum_{i=1}^n \sum_{y,w\in C_i} p_n(y)p^{[n]}(y) p^{[n]}(w)
f(w) \$=
  \sum_{i=1}^n \sum_{y,w\in C_i} p_n(w)p^{[n]}(y) p^{[n]}(w)
f(w) =
  \sum_{i=1}^n \sum_{y\in C_i} p^{[n]}(y) \sum_{w\in C_i} p^{[n+1]}(w)
f(w) =
  E_{n+1}(f)\calcat x.
  \proofend
  $$

We are now ready to  present our main Theorem on ground states.  Unlike
\lcite{\KMSFactors} one cannot prove that all ground states factor
through the conditional expectation $G$ of
\lcite{\MainCondExpectation}.  For example, if we choose the potential
$h$ given by $h_n\equiv1$, then the dynamics is trivial and hence any
state is a ground state, regardless of whether it factors through $G$
or not.  Our result will therefore be restricted to the
characterization of the ground states of the form $\phi\circ G$, where
$\phi$ is a state on $C(X)$.

  \state Theorem 
  {\rm (see \scite{\RenaultThesis}{5.4})}
  \label \MainGround
  Let $X$ be a compact Hausdorff space and $\R=\{R_n\}_{n\in\N}$ an
approximately proper equivalence relation on $X$ such that the
quotient map relative to each $R_n$ is a covering map.  Fix a strictly
positive potential $p = \{p_n\}_{n\in\N}$ satisfying
\lcite{\ConditionOnPForCondexp} and let $E_n$ be the conditional
expectations provided by \lcite{\ConditionOnPForCondexp}.  Also let
$\sigma$ be a one-parameter group of automorphisms of $\Alg$ obtained
from a strictly positive potential $h$.  Given a measure $\mu$ on $X$ let
$\phi$ be the state on $C(X)$ given by integration against $\mu$.  Then
the composition $\psi = \phi\circ G$ is a ground state on $\Alg$ if and only if
the support of $\mu$ is contained in $M^\R_h$.

  \proof
  Let $a,b,c,d\in C(X)$, let $n,m\in\N$, and let $z=\a+i\b$.  If
$n\leq m$ we have by \lcite{\AlgFormula} that
  $$
  \psi\big((a e_n b) \sigma_z(c e_m d)\big) =
  \psi\( a E_n\(b c h^{i \a[n]} h^{-\b[n]} \)e_m h^{-i \a[n]}h^{\b[n]}
d \) \$=
  \int a E_n\(b c h^{i \a[n]} h^{-\b[n]} \)\lambda^{-[m]} h^{-i
\a[n]}h^{\b[n]} d \d \mu =
  \int f E_n^\b(g) \d \mu,
  \eqno{(\dagger)}
  $$
  where $f = a \lambda^{-[m]} h^{-i \a[n]} d$, 
  $g=b c h^{i \a[n]}$,
  and $E_n^\b$ is defined as in \lcite{\TechnicalGround} in terms of
$h^{[n]}$.

If $n\geq m$ we instead have 
  $$
  \psi\big((a e_n b) \sigma_z(c e_m d)\big) =
  \int f E_m^\b(g) \d \mu,
  \eqno{(\ddagger)}
  $$
  where $g$ is as above and $f$ is now $a \lambda^{-[n]} 
h^{-i \a[n]}d$.

Assuming that the support of
$\mu$ is contained  in $M^\R_h$ it follows from
\lcite{\TechnicalGround.ii} that both $(\dagger)$ and $(\ddagger)$ are
bounded as $z$ runs in the upper half plane and hence that $\psi$ is a
ground state.  The converse  also follows easily from
\lcite{\TechnicalGround}.
  \proofend


  \section {Ground states and maximizing measures}
  Consider a fixed Holder real function $H>0$.  We say $\tilde{H}$ is
cohomologous to $H$ if there exists a real function $V$ and real
constant $c$ such that $\log \tilde{H}= \log H - [(V \circ T) - V]
+c$.

An important point in section \lcite{\ThermoSection} is that for a
given $\beta$ the measure $\nu_{H,\beta}$ is an eigenmeasure and
therefore not necessarily invariant for the expanding transformation
$T$.  Given $H$ there exists however another potential $\tilde{H}$,
cohomologous to $H$ such that the eigenmeasure $\nu_{\tilde{H},\beta}$
is an invariant measure.

We would like to investigate similar properties for the ground state
problem. In principle, it can happen that for a certain $H$ there is
no invariant measure $\mu$ with support inside $M_H$ of Theorem
\lcite{\MainGround}.

Given $H$ it will follow from our reasoning in this section that the
measure $\mu$ of Theorem \lcite{\MainGround} associated to a certain
$\tilde{H}$ (cohomologous to $H$) is a maximizing measure in the sense
of \cite{\CLT} and therefore invariant.  These measures are the
analogous (for the case of expanding maps) of the Aubry-Mather
measures of Lagrangian Mechanics.  In the case of the geodesic flow in
compact surfaces of negative curvature they exactly correspond under
the action of the discrete group of Moebius tranformations in the
boundary of the Poincare disk (see \cite{\BS} and \cite{\LT}).

We denote by $ {\cal M} (T)$ the set of invariant probabilities for $T$.
 
First we will recall some general results for maximizing measures.
 
  \definition Given an $\alpha$-Holder function $B$ we denote $$
\hbox{Hol}_\alpha (B)= \sup_{d(x,y) >0} \left\{ \frac{| B(x) -
B(y)|}{d(x,y)^\alpha} \right\}.$$ If we denote by $|| B||_\infty$ the uniform norm,
then we define the $\alpha$-Holder norm of $B$ by $||B ||_\alpha=
\hbox{ Hol}_{\alpha} (B) + || B||_\infty.$  We also let ${\cal
H}_\alpha$ be the set of $\alpha$-Holder functions.

  \definition Given $\log H\in {\cal H}_\alpha$ we define $$ m(H)=\sup
\left\{ - \int \log H(x) d \rho(x) \, |\, \rho \in {\cal M} (T) \right\} $$ and
$${\cal M}_H (T)= \left\{\rho \in {\cal M} (T) : -\int \log H(x) d
\rho (x) = m(H) \right\}.$$
  We call any $\rho\in
{\cal M}_H (T)$ a maximizing measure for $H$ and it will generically
denoted by $\mu_H$.

It is shown in \scite{\CLT}{Proposition 15} that a measure $\mu$ is
maximizing if and only if its support is contained in the
$\Omega(-\log H,T)$ set (see \cite{\CLT} for definition). This result
is the version of Theorem \lcite{\MainGround} above for the case of
invariant measures.
We refer the reader to \cite{\CLT} for general references on the topics
considered in the present section.

Consider ${\cal F}_\alpha^{+}= \cup_{\gamma > \alpha} {\cal H}_\gamma$ 
equiped with
the $\alpha$-norm.

  \state Theorem 
  \label \DozePontoUm
  {\rm (\scite{\CLT}{page 1382})} For an open and dense set 
${\cal G}$ contained in ${\cal F}_\alpha^{+}$,
when $-\log H\in {\cal G}$ then the measure $\mu_H\in {\cal M}_H (T)$ 
is unique and has support in an unique periodic orbit.

It can be shown that for any $H$, the omega-limit set 
of  points in  $M_H$ (of Theorem \lcite{\MainGround}) is contained in the
support of the maximizing measure $\mu_H$. Note that $M_H$ is not
an  invariant
set  for $T$.

In \cite{\CLT} it is shown examples of $H$ where $\mu_H$ is uniquely ergodic and has positive 
entropy.

Assume  $T$ is an expanding transformation on $X$ with 
degree $k$ and  $p=1/k$ as in section \lcite{\ThermoSection}.
We will consider the associated $C^{*}$-algebra $\Alg$ as before. 

Suppose $H$ is strictly positive and Holder and consider the
corresponding $\sigma_t$.

We will say that a measure $\nu$ is a \stress{ground measure} when the
state on $\Alg$ given by $\phi=\nu\circ G$ is a ground state, as in
Theorem \lcite{\MainGround}.

Note that a measure is maximizing for  $H$ Holder, 
if and only if, it is maximizing for $\tilde{H}$ Holder, where 
$-\log \tilde{H}$
cohomologous to $-\log H$. We will also describe the measure $\nu$ associated 
to $\tilde{H}$ as a maximizing measure for $- \log H$ 
(or for $-\log \tilde{H}$)  in the sense of \cite{\CLT}.

More precisely, we will show that one can find $V>0$, Holder such that 
$\tilde{H}(x)=H(x) \,e^{-V(x)+(V\circ T)(x)}$ has a ground
 measure $\nu=\mu_H$, in the sense that, for all 
$f,g\in C(X)$, all $m$ and all complex $\beta$ such that
$\hbox{Re}(\beta)\geq 0$, we have 

$$|\phi (M_{g}\, \sigma_{\beta} (S^{n} (S^{*})^{m} M_{f}))|\leq$$ $$
\int |g \, \, \alpha^m ({\cal L}_p^m (f H^{-\beta [m] } e^{(V
[m]-V\circ T[m]) \beta}) \, )\,\, H^{\beta [m] } e^{(-V [m]+ V\circ T
[m])\beta}| d \nu< || f||_\infty || g||_\infty< \infty,$$ where
$\sigma_z$ is defined by $\tilde{H}$.

We will show that such $\nu$ is invariant and for a generic $H$ it
will follow from Theorem \lcite{\DozePontoUm} that $\nu$ has support
in a unique periodic orbit.

We denote from now on $m(-\log H)= \sup \{ \int- \log H d \rho | \rho
\in {\cal M} (T)\}$.

Note that $m(-\log H+ V- V\circ T)= \sup \{ \int- \log H d \rho | \rho
\in {\cal M} (T)\}= m(- \log H)$, because we are considering $\rho$ an
invariant measure.

By \cite{\CLT}, there exist $V:X \to {\bf R}$, Holder continuous
strictly positive and satisfying for all $x$ the inequality $$V(T(x))
- V(x) \geq - \log H(x) - m(-\log H).$$

This inequality is called a sub-cohomological equation.

The inequality is an equality for $x$ in the support of $\mu_H$.

The function $V$ is defined by
$$V(x) = \sup \{ \sum_{j=0}^{n-1} (-\log H - m (-\log H)) (T^j (y)) \,| 
\,T^n (y) = x, n \in {\bf N}  \}
$$

For $z\in X$, and $n \in {\bf N}$, denote $x_n^i (z)$, $i \in \{1,2,.., 
k^n \},$ the $k^n $ solutions of $T^n (z) =x$.

Fix a point $x$ from now on.

We are going to define a sequence of points $y_n$ inductively.
We set $y_0=x$, and for $y_1$, we choose a point  over the set  
$\{z | \, T(z)=y_0 \}$ such
that $V(T(y_1)) - V(y_1) = -\log H(x) - m (-\log H)$.

From the definition of $V$ one can easily show
that there is always such point $y_1$.

Inductively, given $y_i$, for $y_{i+1}$, we choose a point over the set  
$\{z| T(z)=y_i \}$ such
that $V(T(y_{i+1})) - V(y_{i+1}) = -\log H(y_{i+1}) - m (-\log H)$.

Note that $T(y_{i+1} )= y_i$, for all $i$.

Consider $\mu_n = \frac{1}{n} \sum_{l=0}^{n-1} \delta_{y_{l}} $, and by compactness a 
measure $\nu$ such that is a weak limit $\nu= \lim_{r\to \infty} \mu_{n_r}$.

This is our candidate for being a ground measure for $\tilde{H}=H\,  e^{-V+ V \circ T}$.

We assume from now on  that $H$ is such that $\mu_H$ is unique 
(and uniquely ergodic from \cite{\CLT}).

  \state Proposition $\nu= \mu_H$.

  \proof
$V(x)$ is Holder continuous on $x$, therefore bounded, then 
$$-\int \log H d \nu =-\lim_{r\to \infty} \int \log  H d \mu_{n_r}=-
\lim_{r\to \infty}  \frac{1}{n_r}  \sum_{j=0}^{n_r-1} (\log H (y_j))=  $$
$$\lim_{r\to \infty} \frac{1}{n_r}  \sum_{j=0}^{n_r-1} ( V(y_{j-1}) - V(y_j)+m(-\log H)) =  $$
$$\lim_{r\to \infty}  \frac{1}{n_r}   (V(x)-   V(y_{n_{r}- 1 }) +  
n_r m(-\log H))= m(-\log H).$$

Therefore, $\nu =\mu_H$ and does not depend on $x$. 

We denote such $\nu$ by $\nu_\infty$.
This measure is invariant.

Consider $\phi$ the state satisfying:
for all $m\in {\bf N}$
$$\phi(M_f S^m (S^{*})^m )= \int \frac{f}{\Lambda^{[m]}}  d \nu_\infty.$$

We are interested in $\phi$ such that it is a ground state for $\tilde{ H}=H e^{V- V \circ T}.$

  \state Proposition
  For any Re($\beta)\geq 0$, $m\in {\bf N}$ and $f,g \in C(X)$

$$\int | g \, \,\alpha^m {\cal L}_p^m (f H^{-\beta [m]} e^{(V [m]-
V\circ T [m])\beta} ) H^{\beta [m]} e^{(-V [m]+ V(T)) [m])\beta}| d
\nu_\infty <$$ $$|| g||_\infty|| f||_\infty$$

The proof is similar to Lemma \lcite{\TechnicalGround} using the fact
that from the cohomological equation for any $m$ $$|\frac{
H(y_j)^{\beta [m]} e^{-V (y_j) \beta}}{ H( x_i^m \, (T^m (y_j)) \,
)^{\beta [m]} e^{-V ( x_i^m \, (T^m (y_j))) \beta}} | \leq 1 $$

The conclusion is that
the minimizing measure $\mu_H=\nu_\infty$ determines the ground state $\phi_{\nu_{\infty}}$

It follows from this propostion that

  \state Theorem Given $H>0$ Holder, there is $V>0$ Holder, such that
if $\nu_\infty$ is the maximizing measure for $-\log H$, then the
state $\phi$ defined by $$\phi(M_f S^m (S^{*})^m )= \int
\frac{f}{\Lambda^{[m]}} d \nu_\infty,$$ for all $m\in {\bf N}$, $f\in
C(X)$, is a ground-state for the potential $\tilde{ H}= H e^{-V +
V\circ T}$.

The conclusion is that, if one considers $p=1/k$ and $H>0$, then the
state $\phi_{\nu_\infty}$ associated to Aubry-Mather measure
$\nu_\infty$ for $H$ is a ground-state for some $\tilde{H}$ (such that
$\log \tilde{H}$ is cohomologous to $\log H$).

  \section{Phase transitions}
  We consider in this section an interesting example of KMS state for $\b=1$
associated with the shift $T$ in $2$ symbols $\{0,1\}$ and $p=1/2$.
We will define a special potential $H$. We will show that not always
the equilibrium measures (Statistical Mechanics) for the pressure are
associated to KMS states (Quantum Statistical Mechanics).
  We refer the reader to \cite{\H}, \cite{\Ltwo}, \cite{\Lthree}, \cite{\FL},
\cite{\Y}, \cite{\Lone} for references and results about the topics
discussed in this section.

We are going to introduce the Fisher-Fedenhorf model of Statistical Mechanics
in the therminology of Bernouli spaces and Thermodynamic Formalism \cite{\H}.
  We define  $X$ to be the shift space $X= \{ 0,1\}^\N$
and denote by $T: \ X \to X$  the left shift map. 
We write $z=(z_0 z_1\dots)$ for a point in $X$ 
and $$[w_0 w_1 \dots w_k]= \{ z:\ z_0=w_0, z_1=w_1, \dots, z_k=w_k\}$$ for a
cylinder set of $X$.


For $k>1$ we denote by $M_k \subset X $ the cylinder set
$[\underbrace{111\dots 1}_{k} 0]$ and by $M_0$ the cylinder set $[0].$
The ordered collection $\{M_k\}_{k=0}^\infty$ is a partition of $X$; in
other words these
 sets are disjoint and their
union is the whole space (minus the point $(111\dots)$).
Note that $T$ maps $M_k$ bijectively onto $M_{k-1}$ for $k\geq 1$,
and onto $X$ for $k=0$.
  Also note that 
the point $(111\ldots)$ is fixed by $T$.

For a fixed real constant $\gamma>1$
we consider the   function  $g$ on $X$ such that $g(111\dots .)=0$, 
$$
g(x)=a_k:=-\gamma\log\biggl(\frac{k+1}{k}\biggr),$$
for 
$x\in M_k$, for $k\neq 0,$  and 
$$
g(x)=a_0:=-\log(\zeta(\gamma)),$$ 
for $x\in M_0,$
where $\zeta$ is the Riemann zeta function.


By definition, $$ \zeta(\gamma)=(1^{-\gamma}+2^{-\gamma}+\dots)$$ and
so the reason for defining $a_0$ in such a way is that, if we define
$s_k=a_0 + a_1 +\dots + a_k$, then $\Sigma e^{s_k}=1.$

From now on we assume $\gamma>2$, otherwise we have to consider
sigma-finite measures  and not probabilities in our problem.

The potential $1<(\frac{k+1}{k})^\gamma = H(x)= e^{-g(x)},$ for $x \in
M_k,$ is not H\"older and in fact is not of summable variation. Note
that $H( 111\ldots)=1$. The pressure $P(-\log H) =P(g)=P( \log p +
\log 2- 1\, \log H)=0$ and one can show that there exist two
equilibrium states for such a potential $g$ (in the sense of
minimizing measures for the variational problem): a point mass (the
Dirac delta $\delta{(111\ldots)}$) at $(111\dots )$, and a second
measure which we shall denote by $\tilde{\mu}$ (see \cite{\H}).

The existence of two invariant probabilities: $\tilde{\mu}$ and
$\delta_{(111\ldots)}$; for the variational problem of pressure 
defines what is called a phase transition
in the sense of Statistical Mechanics \cite{\H}, \cite{\Lthree}.

We will describe bellow how to define this measure $\tilde{\mu}$.

Consider as in \cite{\H}
${\cal  L}_g^*$,  
the dual of the 
Ruelle-Perron-Frobenius operator $ {\cal  L}_g$ associated to $g$,
where 
the action of $ {\cal  L}_g$
 on continuous functions is given by
$$  {\cal  L}_{\beta=1} (\phi)(y)=\sum_{T(x)=y}e^{g(x)}\phi(x).$$

The function $P(-\beta \log H) = P(\beta  g)$ is strictly
monotone for $\beta<1$ and constant equal zero for $\beta>1$ (see \cite{\H},\cite{\Lone}).

We claim that there is a unique 
probability measure $\nu$ on $X$ 
which satisfies $ {\cal  L}_g^*\nu=\nu$ \cite{\FL}, \cite{\H}.
To prove this, note first that $\nu$ cannot have any mass at $(111\dots)$;
it follows that $M_0$ has positive mass, and the stipulation that 
$\nu$ be an eigenmeasure then gives a recurrence relation  for the 
masses of $M_k$. Since $T(M_k)=M_{k-1}$ for $k\geq 1$, we have 
that the masses of the sets in this partition are   
 $$\nu(k)= \nu(M_k)=e^{s_k}=\frac{(k+1)^{-
\gamma}}{\zeta(\gamma)},\, k \geq 0 .$$

In particular,
$$\nu(0)=\nu(M_0)=e^{s_0}=e^{a_0}=\frac{1}{\zeta(\gamma)}.$$

By the same reasoning, $\nu$ is determined on all
higher cylinder sets for the partition  $(M_k)_{k=0}^\infty$.
Hence  $\nu $ exists and is unique but not invariant.

The measure $\nu$ defined above is the unique eigenmeasure
for ${\cal  L}_{g}^*$.

The measure defined by the delta-Dirac on $(111\ldots)$ is invariant
but is not a fixed eigenmeasure for ${\cal  L}_g^*$.

There exists $f$ such that ${\cal  L}_g (f)=f$ and
$\log {\tilde H}= (f \circ T - f) - g$
defines  ${\tilde H}$ cohomologous to $e^g$ (see \cite{\H}).

$-\log {\tilde H}$ defines
for the pressure $P(-\log {\tilde H} )$ the same 
two equilibrium measures (as for $g=-\log H$): 
the invariant measure $\tilde{\mu}=f d \nu$ and the delta-Dirac on $(111\ldots)$.
${\cal  L}_H^*$ has a unique eigenmeasure $\tilde{\mu}=f d \nu$ 
which is invariant.

This measure $\tilde{\mu}$ defines a KMS state 
$\phi$ for such $\tilde{H}$, $\beta=1$. 

We can conclude from the above considerations that not always an
equilibrium probability $\rho$ for the pressure is associated to a KMS
state $\phi_\rho$ whithout the hypothesis that $H$ is Holder and
$H>1$. In the present example, this happen because
$\rho=\delta_{(111\ldots)}$ is not an eigenmeasure of ${\cal L}_H^* $
but it is an equilibrium measure for $P(-\beta \log \tilde{H})$ with
$\beta=1$.

In \cite{\Ltwo} and \cite{\Lthree} the lack of differentiability of the Free energy
is analyzed  and in  \cite{\Lthree}, \cite{\FL}, \cite{\Y} it is shown that such systems
present polynomial decay of correlation. In \cite{\Lone}
it is presented a dynamical model with three equilirium states.


  \endgroup

  \bye